\documentclass{amsart}
\usepackage{graphicx}
\usepackage[latin2]{inputenc}
\usepackage{amsmath,amssymb,amsbsy,amsthm,amsfonts}
\usepackage{epsfig} % postscript figures
\begin{document}
%\numberwithin{equation}{section} \marginparwidth=2cm
\def\note#1{\marginpar{\small #1}}

\def\tens#1{\pmb{\mathsf{#1}}}
\def\vec#1{\boldsymbol{#1}}

\def\norm#1{\left|\!\left| #1 \right|\!\right|}
\def\fnorm#1{|\!| #1 |\!|}
\def\abs#1{\left| #1 \right|}
\def\ti{\text{I}}
\def\tii{\text{I\!I}}
\def\tiii{\text{I\!I\!I}}

\def\diver{\mathop{\mathrm{div}}\nolimits}
\def\grad{\mathop{\mathrm{grad}}\nolimits}
\def\Div{\mathop{\mathrm{Div}}\nolimits}
\def\Grad{\mathop{\mathrm{Grad}}\nolimits}

\def\tr{\mathop{\mathrm{tr}}\nolimits}
\def\cof{\mathop{\mathrm{cof}}\nolimits}
\def\det{\mathop{\mathrm{det}}\nolimits}

\def\lin{\mathop{\mathrm{span}}\nolimits}
\def\pr{\noindent \textbf{Proof: }}
\def\pp#1#2{\frac{\partial #1}{\partial #2}}
\def\dd#1#2{\frac{\d #1}{\d #2}}

\def\T{\mathcal{T}}
\def\R{\mathcal{R}}
\def\bx{\vec{x}}
\def\be{\vec{e}}
\def\bef{\vec{f}}
\def\bec{\vec{c}}
\def\bs{\vec{s}}
\def\ba{\vec{a}}
\def\bn{\vec{n}}
\def\bphi{\vec{\varphi}}
\def\btau{\vec{\tau}}
\def\bc{\vec{c}}
\def\bg{\vec{g}}
\def\bd{\vec{d}}

\def\bW{\tens{W}}
\def\bT{\tens{T}}
\def\bD{\tens{D}}
\def\bF{\tens{F}}
\def\bB{\tens{B}}
\def\bV{\tens{V}}
\def\bS{\tens{S}}
\def\bI{\tens{I}}
\def\bi{\vec{i}}
\def\bv{\vec{v}}
\def\bfi{\vec{\varphi}}
\def\bk{\vec{k}}
\def\bh{\vec{h}}
\def\b0{\vec{0}}
\def\bom{\vec{\omega}}
\def\bw{\vec{w}}
\def\p{\pi}
\def\bu{\vec{u}}

\def\ID{\mathcal{I}_{\bD}}
\def\IP{\mathcal{I}_{p}}
\def\Pn{(\mathcal{P})}
\def\Pe{(\mathcal{P}^{\eta})}
\def\Pee{(\mathcal{P}^{\varepsilon, \eta})}

\def\Ln#1{L^{#1}_{\bn}}

\def\Wn#1{W^{1,#1}_{\bn}}

\def\Lnd#1{L^{#1}_{\bn, \diver}}

\def\Wnd#1{W^{1,#1}_{\bn, \diver}}

\def\Wndm#1{W^{-1,#1}_{\bn, \diver}}

\def\Wnm#1{W^{-1,#1}_{\bn}}

\def\Lb#1{L^{#1}(\partial \Omega)}

\def\Lnt#1{L^{#1}_{\bn, \btau}}

\def\Wnt#1{W^{1,#1}_{\bn, \btau}}

\def\Lnd#1{L^{#1}_{\bn, \btau, \diver}}

\def\Wntd#1{W^{1,#1}_{\bn, \btau, \diver}}

\def\Wntdm#1{W^{-1,#1}_{\bn,\btau, \diver}}

\def\Wntm#1{W^{-1,#1}_{\bn, \btau}}

%------------------------------------------------
\newtheorem{Theorem}{Theorem}[section]
%{\theorembodyfont{\rmfamily} \newtheorem{Example}{Example}}
\newtheorem{Example}{Example}[section]
\newtheorem{Lem}{Lemma}[section]
\newtheorem{Rem}{Remark}[section]
\newtheorem{Def}{Definition}[section]
\newtheorem{Col}{Corollary}[section]

%\tableofcontents

\numberwithin{equation}{section}

\title[On Kolmogorov's two-equation model of turbulence]{Large data analysis for Kolmogorov's two-equation model of turbulence}\thanks{Miroslav Bul\'{\i}\v{c}ek and Josef M\'{a}lek research is supported by the Ministry of Education, Youth and Sport through the ERC-CZ project LL1202.}

\author[M. Bul\'{\i}\v{c}ek]{Miroslav Bul\'{\i}\v{c}ek}
\address{Charles University, Faculty of Mathematics and
Physics, Mathematical Institute\\ Sokolovsk\'{a}~83,
186~75~Prague~8, Czech~Republic}
\email{mbul8060@karlin.mff.cuni.cz}

\author[J. M\'{a}lek]{Josef M\'{a}lek}
\address{Charles University, Faculty of Mathematics and
Physics, Mathematical Institute\\ Sokolovsk\'{a}~83,
186~75~Prague~8, Czech~Republic}
\email{malek@karlin.mff.cuni.cz}

\keywords{two-equation model of turbulence, k-epsilon model, existence, weak solution, suitable weak solution, Navier's slip, degenerate viscosity}
\subjclass[2000]{35Q30, 35Q35, 76F60}

\begin{abstract}
Kolmogorov seems to have been the first to recognize that a two-equation model of turbulence might be used as the basis of turbulent flow prediction. Nowadays, a whole hierarchy of phenomenological two-equation models of turbulence is in place. The structure of their governing equations is similar to the Navier-Stokes equations for incompressible fluids, the difference is that the viscosity is not constant but depends on the fraction of the scalar quantities that measure the effect of turbulence: the average of the kinetic energy of velocity fluctuations (i.e. the turbulent energy) and the measure related to the length scales of turbulence. For these two scalar quantities two additional evolutionary convection-diffusion equations are augmented to the generalized Navier-Stokes system. Although Kolmogorov's model has so far been almost unnoticed it exhibits interesting features. First of all, in contrast to other two-equation models of turbulence there is no source term in the equation for the frequency. Consequently, nonhomogeneous Dirichlet boundary conditions for the quantities measuring the effect of turbulence are assigned to a part of the boundary. Second, the structure of the governing equations is such that one can find an ``equivalent" reformulation of the equation for turbulent energy  that eliminates the presence of the energy dissipation acting as the source in the original equation for turbulent energy and which is merely an $L^1$ quantity. Third, the material coefficients such as the viscosity and turbulent diffusivities may degenerate, and thus the a priori control of the derivatives of the quantities involved is unclear.

We establish long-time and large-data existence of a suitable weak solution to three-dimensional internal unsteady flows described by Kolmogorov's two-equation model of turbulence. The governing system of equations is completed by initial and boundary conditions; concerning the velocity we consider generalized stick-slip boundary conditions. The fact that the admissible class of boundary conditions includes various types of slipping mechanisms on the boundary makes the result robust from the point of view of possible applications.
\end{abstract}

\maketitle

\section{Introduction}

In this paper, we establish long-time and large data existence of suitable weak solution to an initial and boundary-value problem associated to a nonlinear system of PDEs proposed, in 1942, by A.~N.~Kolmogorov to describe three-dimensional unsteady turbulent flows, see \cite{Kolmogorov1942}\footnote{Two English translations are available, see Appendix in \cite{Spalding1991} and the paper No. 48 in \cite{tichomirov91}}.

We first formulate the problem and provide its reformulation that is equivalent in the context of regular enough functions but exhibits better mathematical properties in the context of weak solutions. We then formulate the assumptions on data following the aim to make them general enough to include relevant physical situations, and state the main result. Next, after a brief introduction to the reduction of Kolmogorov's system to a one-equation model of turbulence, we add our motivation for investigating Kolmogorov's 1942 system of equations. Then, we highlight the main novelties of our result and we conclude this introductory section by recalling relevant mathematical results. In Section \ref{scheme}, we outline the scheme of the proof of the main result and introduce a hierarchy of three levels of approximate problems and formulate the lemmas concerning the existence of weak solutions to these approximate problems. The proofs of these lemmas are given in Section \ref{aupr}. In the final section \ref{SMT} we study the limit of the highest level approximate problem to the original problem, and thus complete the proof of the main theorem.

\subsection{Formulation of the problem} \label{formulation}

Let  $\Omega \subset \mathbb{R}^3$ be an open bounded set with Lipschitz boundary $\partial \Omega$ and let $T>0$ denote the length of the time interval. We set $Q:=(0,T)\times \Omega$. Our goal is to analyze the following problem: find  $(\bv,p, \omega, b):Q \to \mathbb{R}^3\times \mathbb{R}\times  \mathbb{R}_+\times \mathbb{R}_+$ solving Kolmogorov's two-equation model of turbulence (see \cite{Kolmogorov1942}) that takes the following form:
\begin{align}
\diver \bv &=0, \label{BM}\\
\partial_{t} \bv + \diver (\bv \otimes \bv) - 2 \nu_0\diver \left( \frac{b}{\omega}\bD(\bv) \right) &= -\nabla p,\label{BLM}\\
\partial_{t} \omega + \diver (\omega\bv) - \kappa_1\diver \left (\frac{b}{\omega} \nabla \omega \right) &= -\kappa_2 \omega^2,\label{len}\\
\partial_{t} b +\diver (b\bv)-\kappa_3\diver \left(\frac{b}{\omega} \nabla b\right)&=-b\omega + \kappa_4 \frac{b}{\omega}|\bD(\bv)|^2.\label{TKE}
\end{align}
Here, $\bv$ stands for the average velocity of the fluid\footnote{Let $(\tilde{\bv},\tilde p)=(\tilde{\bv}(t,x), \tilde{p}(t,x))$ denote the solution to the Navier--Stokes equations. Let further $\bv$ denote a (time, spatial or stochastic) average of $\tilde{\bv}$, , i.e. $\bv:= \langle \tilde{\bv} \rangle$ where brackets $\langle \cdot \rangle$ denotes here a certain averaging. Then $\tilde{\bv} = \bv + \bv'$, where $\bv'$ denotes the velocity of fluctuations. Denoting $b:= \tfrac13 \langle |\bv'|^2  \rangle$ one observes that $b$ is related to the average of the turbulent kinetic energy defined through $b:= \tfrac13 \langle |\bv'|^2  \rangle$ via the equation $b= \tfrac23 k$. Finally the average pressure $p$ is set to be $p:= \frac{\langle \tilde{p} \rangle}{\varrho_*} + b$.} and $\bD(\bv)$ is the symmetric part of its gradient, $b$ denotes $\tfrac{3}{2}$ of the turbulent kinetic energy (i.e., the kinetic energy of the velocity fluctuations), $p$ is the sum of $b$ and the average of the mean normal stress divided by the constant density and $\omega$ is the frequency related to the length scale $\ell$ by the relation $\omega:=c\sqrt{b}/\ell$, where $c>0$ is a constant. For simplicity we neglect external body forces. In \eqref{BM}--\eqref{TKE} the material parameters  $\nu_0,\kappa_1, \ldots, \kappa_4$ are assumed to be  given positive constants; Kolmogorov specified $\kappa_2$ to be $\tfrac{7}{11}$ and considered $\kappa_4 = 2\nu_0$, see \cite{Kolmogorov1942}.

To complete the system \eqref{BM}--\eqref{TKE} we need to specify the initial and boundary data. Regarding the initial conditions, we assume that
\begin{equation}
\bv(0,x)=\bv_0(x), \qquad b(0,x)=b_0(x), \qquad \omega(0,x)=\omega_0(x) \qquad \textrm{ for } x\in \Omega; \label{ID}
\end{equation}
we shall put further restrictions on the given $\bv_0$, $b_0$ and $\omega_0$ in the next subsection.

Concerning the boundary conditions for the turbulent kinetic energy $b$ and the frequency $\omega$ we first notice  that there is no source term in the equation $\omega$ which indicates (see Spalding \cite{Spalding1991} and Subsect. \ref{history} below) that these turbulent quantities have to be generated on some part of the boundary. This is why
we assume that $\partial \Omega$ consists of two smooth open disjoint parts $\Gamma$ and $\Gamma^c$ such that $\overline{\Gamma \cup \Gamma^c}=\partial \Omega$ and we
consider mixed boundary conditions of the form
\begin{align}
b &=b_{\Gamma} \quad &&\textrm{on } (0,T)\times \Gamma,\label{bc3}\\
\frac{b}{\omega}\nabla b \cdot \bn&=0 \quad &&\textrm{on } (0,T)\times \Gamma^c, \label{bc4}
\end{align}
and
\begin{align}
\omega &=\omega_{\Gamma} \quad &&\textrm{on } (0,T)\times \Gamma,\label{bc5}\\
\frac{b}{\omega}\nabla \omega \cdot \bn&=0 \quad &&\textrm{on } (0,T)\times \Gamma^c. \label{bc6}
\end{align}

Next, we focus on the boundary conditions for the velocity field. We will consider internal flows, i.e., we assume that
\begin{align}
\bv \cdot \bn &=0 \quad &&\textrm{on } (0,T)\times \partial \Omega,\label{bc1}
\end{align}
where $\bn$ denotes the unit outer normal vector on $\partial \Omega$. Let further, for any vector $\bw$ originating at the point $x\in \partial \Omega$, $\bw_{\btau}:= \bw - (\bw\cdot \bn(x))\bn(x)$ denote the projection of $\bw$ on the tangent plane of $\partial \Omega$ at $x$. Using the notation
$$
\bs:= - 2\nu_0 \left(\frac{b}{\omega}\bD(\bv) \bn\right)_{\btau}
$$
for the projection of the normal traction\footnote{Note that $\bs$ also equals to $(\bT\bn)_{\btau}$ where $\bT = -p\bI + 2\nu_0 \frac{b}{\omega} \bD(\bv)$ is the Cauchy stress tensor associated with \eqref{BM}--\eqref{TKE}.}, we can formulate the boundary condition relating $\bs$ to $\bv_{\btau}$. Note that the precise form of this boundary condition is a subject of intense investigations, particularly for turbulent flow, see \cite{ChaLe14}. We aim at including various slipping mechanisms as well as an activated transition from no-slip to partial slip where the threshold itself can depend on the kinetic turbulent energy, on the mixing length and also on the spatial and time variable to cover the case of different material properties on solid boundaries. Thus our set of assumptions on $\bs$ is given through the following condition:
\begin{equation}
\left.
\begin{aligned}
|\bs(t,x)|&\le \sigma(t,x,b,\omega) &&\implies \bv_{\btau}(t,x) = \b0,\\
|\bs(t,x)|& > \sigma(t,x,b,\omega) &&\implies \bs(t,x) = \bg (t,x,b,\omega, \bv_{\btau})
\end{aligned}
\right\} \textrm{ on } (0,T)\times \partial \Omega, \label{bc2}
\end{equation}
where $\sigma: (0,T)\times \partial \Omega \times \mathbb{R}_+^2 \to \mathbb{R}_+$ is a given (continuous) threshold function and $\bg:(0,T)\times \partial \Omega \times \mathbb{R}_+^2 \times (\mathbb{R}^3 \setminus \{\b0\})\to \mathbb{R}^3$ denotes the possible (continuous) slip function, which may not be defined for $\bv_{\btau} = \b0$ in order to be able to cover general threshold slip.

Note that when $\sigma \to \infty$, the condition \eqref{bc2} approximates the no-slip boundary condition
\begin{equation*}
  \bv_{\btau}= \b0.
\end{equation*}
On the other hand, setting $\sigma=0$, \eqref{bc2} includes as a special case Navier's slip boundary condition described by
\begin{equation*}
\bs = \gamma_* \bv_{\btau},% \quad \textrm{ on } (0,T)\times \partial \Omega,
\end{equation*}
where $\gamma_*>0$ is the friction coefficient. Finally, \eqref{bc2} also includes the standard stick-slip boundary condition
\begin{equation*}
\begin{split}
|\bs| &\le \sigma_{*} \Leftrightarrow \bv_{\btau} = \b0, \\
|\bs| &>\sigma_{*} \Leftrightarrow \bs = \sigma_{*} \frac{\bv_{\btau}}{|\bv_{\btau}|} + \gamma_* \bv_{\btau},
\end{split}
\end{equation*}
where the threshold $\sigma_*$ is a positive constant and the fluid slips along the boundary as in the case of Navier's slip boundary condition once the tangent projection of the normal traction $\bs$ exceeds the threshold $\sigma_*$.

We can think of \eqref{bc2} as a (continuous) curve defined on the Cartesian product $\mathbb{R}^3\times \mathbb{R}^3$ in the variables $\bs$ and $\bv_{\btau}$ parametrized by $t$, $x$, $b$ and $\omega$. This means that we can rewrite \eqref{bc2} in the equivalent implicit form
\begin{equation}
\bh (\cdot, b, \omega; \bs, \bv_{\btau}) = \b0 \qquad \textrm{ on } (0,T)\times \partial \Omega. \label{bc2impl}
\end{equation}
In a particular case, when $\bg$ is of the form
\begin{equation*}
\bg (\cdot, b, \omega; \bs, \bv_{\btau}) = \sigma(\cdot, b, \omega) \frac{\bv_{\btau}}{|\bv_{\btau}|} + \gamma(\cdot, b, \omega; \bs, \bv_{\btau}) \bv_{\btau},
\end{equation*}
where $\gamma$ stands for the generalized friction function, the function $\bh$ introduced in \eqref{bc2impl} takes the form
\begin{equation*}
\bh (\cdot, b, \omega; \bs, \bv_{\btau}) = \gamma(\cdot, b, \omega; \bs, \bv_{\btau}) \bv_{\btau} - \frac{\left(|\bs|- \sigma(\cdot, b, \omega)\right)^{+}}{|\bs|}\bs,
\end{equation*}
where $\cdot$ stands for $t,x$ and $z^{+}:= \max\{0,z\}$.

We finally note that we could consider threshold conditions for $b$ and $\omega$ (similar to \eqref{bc2} for the velocity field), which could be even more appropriate
on the part $\Gamma^c$ of the boundary in order to describe activated occurance of turbulence. Another direction worthy of investigation is to to consider the dependence of $b_{\Gamma}$ and $\omega_{\Gamma}$ on the velocity field, see \cite[Section 5]{ChaLe14}. We however do not study these two generalizations here.
We also note that the mathematical theory regarding the generalized stick-slip conditions \eqref{bc2} developed in this paper extends recent mathematical approaches \cite{bulcek-malek-a:2016, bulcek-malek-b:2016} involving stick-slip boundary conditions. In \cite{bulcek-malek-a:2016,bulcek-malek-b:2016}, we require that $\bg$ depends merely on $\bv_{\btau}$ and in addition this dependence is monotone.

\subsection{Difficulties and the ``equivalent" reformulations of \eqref{len} and \eqref{TKE}}

The system \eqref{BM}-\eqref{TKE} consists of the generalized Navier-Stokes equations coupled with two scalar evolutionary convection-diffusion equations. The quantity $\frac{b}{\omega}$ measures the effective kinematic viscosity and the effective diffusivity of turbulence. It seems reasonable to assume that the initial and boundary data for the frequency $\omega$ are uniformly positive and bounded from above, which together with the structure of the equation \eqref{len} implies that $\omega$ remains uniformly positive and bounded from above over the whole time cylinder $Q$. On the other hand, $b$ is required to be merely positive initially (which together with the structure of the equation for $b$ implies that $b$ is at least nonnegative in $Q$). Consequently, $\frac{b}{\omega}$ might degenerate\footnote{This feature somehow puts the system \eqref{BM}-\eqref{BLM} between the Navier-Stokes and Euler equations.} and it is not a priori evident that one can control spatial derivatives of $\bv$, $b$ and $\omega$. Note that the formal identity (valid for any $t\in (0,T)$)
\begin{equation}
\int_{\Omega} |\bv(t,x)|^2 \, dx + \int_0^t\int_{\Omega} \frac{4\nu_0 b}{\omega}|\bD(\bv)|^2 \, dx\, ds + 2 \int_0^t\int_{\partial \Omega} \bs \cdot \bv_{\btau} \, dS\, ds = \int_{\Omega}|\bv_0(x)|^2 \, dx, \label{pepa1}
\end{equation}
obtained after integrating, over $(0,t) \times \Omega$,  the equation
\begin{equation}
\partial_{t} |\bv|^2 + \diver (|\bv|^2 \bv) - \diver \left( \frac{4\nu_0 b}{\omega}\bD(\bv)\bv \right) + \frac{4\nu_0 b}{\omega} |\bD(\bv)|^2 = -2\diver (p \bv),\label{pepa2}
\end{equation}
does not imply that $\bD(\bv)$ belongs to $L^2(Q)$. On the other hand, \eqref{pepa1} implies that the last term at the right-hand side of \eqref{TKE} belongs to $L^1(Q)$. Consequently,
\begin{equation}
\sup_{t\in (0,T)} \int_{\Omega} b(t,x) \, dx < + \infty. \label{pepa3}
\end{equation}

Furthermore, multiplying \eqref{TKE} first by $\frac{1}{1+b}$ and then by $\frac{1}{(1+b)^{\varepsilon}}$ for $\varepsilon\in (0,1)$ we obtain, roughly speaking\footnote{One needs to take into account the Dirichlet boundary conditions on $(0,T) \times \Gamma$ as well as possible degeneracy of $b$ in $Q$. The details are give below in the proof of the main result.}, the following estimates:
(valid for any $t\in (0,T)$)
\begin{equation}
\int_{\Omega} |\ln b(t,x)| \, dx + \int_0^t\int_{\Omega} |\bD(\bv)|^2 \, dx\, ds + \int_0^t\int_{\Omega} |\nabla b^{1-\frac{\varepsilon}{2}}|^2 \, dx\, ds < \infty, \label{pepa4}
\end{equation}
and we conclude that $\nabla \bv$ is $L^2$-integrable over $Q$ and $\nabla b$ is almost $L^2$-integrable over $Q$. This result regarding $\nabla \bv$ puts the problem in the same function space setting as the Navier-Stokes equations. However, in contrast to the Navier-Stokes equation, the system \eqref{BM}--\eqref{TKE} requires one to handle much more severe nonlinearities, in particular the last term in \eqref{TKE}, which belongs merely to $L^1(Q)$. Furthermore, although
\begin{equation}
\int_0^t\int_{\Omega} b |\nabla w|^2 \, dx\, ds < \infty \label{pepa5}
\end{equation}
we do not obtain control of $\nabla \omega$. Combining \eqref{pepa4} and \eqref{pepa5} we however control $\nabla (bw)$. This brings us to the conclusion to reformulate \eqref{len} and \eqref{TKE}.

Instead of \eqref{len} we consider
\begin{equation}
\partial_{t} \omega + \diver (\omega\bv) - \kappa_1\diver \left(\frac{\nabla (b\omega) - \omega\nabla b}{\omega} \right) = -\kappa_2 \omega^2.\label{len-pepa}
\end{equation}

The presence of the $L^1(Q)$-nonlinearity in \eqref{TKE} is overcome by replacing \eqref{TKE} by the equation for the sum of the kinetic energy of $\bv$  and the kinetic energy of the velocity fluctuations, i.e., a multiple of $b$. More precisely, setting
\begin{equation}
E:=\frac12 |\bv|^2 +\frac{2 \nu_0}{\kappa_4}b,\label{DFE}
\end{equation}
multiplying \eqref{TKE} by $\frac{\nu_0}{\kappa_4}$, taking the scalar product of $\bv$ and \eqref{BLM}, and finally summing the resulting identities, one arrives at
\begin{equation}
\partial_{t} E + \diver \left( \bv(E+p)\right) - 2\nu_0\diver \left( \frac{\kappa_3 b}{\kappa_4 \omega} \nabla b + \frac{b}{\omega}\bD(\bv)\bv\right) +\frac{2\nu_0}{\kappa_4}b\omega = 0. \label{BKE}
\end{equation}
It is advantageous that all terms in \eqref{BKE} are in the divergence form. On the other hand, the pressure $p$ appears in \eqref{BKE} and its (sufficient) integrability is required in order to handle the convergence in the relevant nonlinear term. Within the setting considered here, this eliminates the no-slip boundary condition from further consideration, see \cite{frehsemalek2002} for more details.

To conclude, within the context of regular enough solution, the system \eqref{BM}--\eqref{TKE} is equivalent to the system consisting of \eqref{BM},
\eqref{BLM}, \eqref{len-pepa} and \eqref{BKE}. Within the context of weak solutions, the latter system has better features and the existence of weak solution to this system will be established in this study.

Moreover, if one requires that a weak solution to \eqref{BM}, \eqref{BLM}, \eqref{len-pepa} and \eqref{BKE} in addition satisfies
\begin{equation}
\partial_{t} b +\diver (b\bv)-\kappa_3\diver \left(\frac{b}{\omega} \nabla b\right) \ge -b\omega + \kappa_4 \frac{b}{\omega}|\bD(\bv)|^2,\label{TKI}
\end{equation}
in a weak sense, then it is natural to call such a solution \emph{suitable weak solution}. Indeed, subtracting \eqref{TKI} from \eqref{BKE}, one deduces that
\begin{equation}
\partial_{t} |\bv|^2 + \diver \left((|\bv|^2 + 2p)\bv \right) - \diver \left(4\nu_0 \frac{b}{\omega} \bD(\bv)\bv\right) + 4\nu_0 \frac{b}{\omega} |\bD(\bv)|^2 \le 0, \label{SWS}
\end{equation}
which is the usual notion of suitable weak solution to the Navier-Stokes system, noticing that if $\frac{b}{\omega}=1$, then
\begin{equation*}
- 4 \diver \left(\bD(\bv)\bv\right) + 4 |\bD(\bv)|^2 = - \Delta |\bv|^2 + 2|\nabla \bv|^2.
\end{equation*}

\subsection{Notation}
We use the standard notation for Lebesgue, Sobolev and Bochner spaces. In order to distinguish between scalar-, vector-, and tensor-valued functions, we use small letters for scalars, small bold letters for vectors and capital bold letters for tensors. Moreover, to simplify the notation for any Banach space $X$ we use the  abbreviation $X^k:= \underset{k - times}{\underbrace{X \times \cdots \times X}}.$ Next, since we need to deal with vector-valued functions having zero normal part on the boundary, we require that $\Omega$ is a Lipschitz domain\footnote{The trace operator is well defined for Lipschitz domains.} and we denote
\begin{align*}
W^{1,r}_{\bn}&:=\left\{ \bv \in W^{1,r}(\Omega)^3: \; \bv \cdot \bn =0 \textrm{ on } \partial \Omega \right\},\\
W^{1,r}_{\bn,\diver}&:=\left\{ \bv \in W^{1,r}_{\bn}: \; \diver \bv =0 \textrm{ in }  \Omega \right\},\\
W^{-1,r'}_{\bn}&:=\left(W^{1,r}_{\bn} \right)^{*}, \quad W^{-1,r'}_{\bn,\diver}:=\left(W^{1,r}_{\bn,\diver} \right)^{*},\\
L^2_{\bn,\diver}&:= \overline{W^{1,2}_{\bn,\diver}}^{\|\, \|_{2}}.
\end{align*}
All of the above spaces are the Banach spaces, which are for $r\in [1,\infty)$ separable and for $r\in (1,\infty)$ reflexive. Next, in order to incorporate Sobolev functions vanishing on a part of the boundary, we denote for an arbitrary smooth relatively open $\Gamma \subset \partial \Omega$
$$
W^{1,r}_{\Gamma}(\Omega):= \{u\in W^{1,r}(\Omega); \; u=0 \textrm{ on } \Gamma \}.
$$
We shall also employ the following notation for functions having zero mean value:
$$
L^r_0(\Omega):=\{u\in L^r(\Omega): \, \int_{\Omega} u \; dx =0\}.
$$
Since we shall also work with sequences that are pre-compact only in the space of measures (bounded sequences in $L^1$),  we denote the space of Radon measures on a set $V$ by $\mathcal{M}(V)$. We will also use the standard notation for dual spaces to spaces of Sobolev functions, i.e., we set $W^{-1,p'}_{\Gamma}(\Omega):=(W^{1,p}_{\Gamma}(\Omega))^*$. Similarly, we denote in the corresponding way also dual spaces to $W^{1,p}_{\bn}$ and $W^{1,p}_{\bn,\diver}$. Finally, in order to simplify the notation, we define the natural energy set for $b$ as follows:
\begin{equation}
\begin{split}
\mathcal{E}:=&\Big\{ b\in L^{\infty}(I;L^1(\Omega)); \; b> 0 \textrm{ a.e. in } Q, \; \ln b \in L^{\infty}(0,T; L^1(\Omega)),\\
&\qquad (1+b)^\lambda \in L^2(I; W^{1,2}(\Omega)) \textrm{ for all } \lambda \in (0,1),\\
&\qquad b-b_{\Gamma} \in L^1(0,T;W^{1,1}_{\Gamma}(\Omega))\Big \}.
\end{split} \label{eee}
\end{equation}
In addition, to shorten the formula we also use the abbreviation $(a,b)_V:=\int_V ab$ whenever $a\in L^r(V)$ and $b\in L^{r'}(V)$ and in particular  if $V=\Omega$ we shall omit writing this subscript in what follows. Similarly, we use the same notation for vector- and tensor-valued functions. In the case of dualities, we will frequently use the abbreviated notation $\langle a, b \rangle:=\langle a, b\rangle_{X,X^*}$ whenever $a\in X$ and $b\in X^*$ and the meaning of the duality pairing is clear from the context.

\subsection{Assumptions on the data}
In this subsection, we specify our requirements on the data. In particular, our goal is to cover the natural case (with the only assumption of bounded energy) and we also want to include the possibility that the turbulent kinetic energy is not uniformly positive initially. In addition, we keep the conditions on the boundary data as general as possible in order to include very general behavior of the Cauchy stress and the velocity field on the boundary.

Thus, we first specify the requirements on the initial data. For the velocity $\bv$ and the turbulent energy $k$ we assume that
\begin{align}
&\bv_0 \in L^2_{\bn, \diver},\label{Assv}\\
&b_0 \in L^1(\Omega), \; b_0> 0 \textrm{ a.e. in } \Omega,\; \ln b_0 \in L^1(\Omega). \label{Assk}\\
\intertext{Next, for the frequency $\omega$, we assume that there exist $0<\omega_{\min}\le \omega_{\max}<\infty$ such that}
&\omega_0 \in L^\infty(\Omega) \quad \textrm{ and } \quad \omega_{\min}\le \omega_0\le \omega_{\max} \qquad\textrm{ a.e. in } \Omega. \label{Asso}
\end{align}
Concerning the boundary conditions for $b$ and $\omega$, we simplify the situation by assuming that $\omega_{\Gamma}$ and $b_{\Gamma}$ can be extended  onto the whole of $Q$ (and we denote these extensions again by $\omega_{\Gamma}$ and $b_{\Gamma}$) such that
\begin{align}
\omega_{\Gamma} &\in L^{\beta_{\Gamma}}(0,T;W^{1,\beta_{\Gamma}}(\Omega)) \cap W^{1,1}(0,T; L^1(\Omega))\qquad \textrm{ for some } \beta_{\Gamma}> \frac{16}{5}, \label{bd1}\\
b_{\Gamma} &\in   L^2(0,T;W^{1,2}(\Omega)) \cap W^{1,1}(0,T; L^1(\Omega)). \label{bd2}
\end{align}
In addition, we require that, for $\omega_{\min}$, $\omega_{\max}$ introduced above and for some $0<b_{\min}\le b_{\max} <\infty$,
\begin{align}
\omega_{\min}&\le \omega_{\Gamma} \le \omega_{\max} &&\textrm{a.e.  in  } Q,\label{bd3}\\
b_{\min}&\le b_{\Gamma}\le b_{\max} &&\textrm{a.e. in  }Q.\label{bd4}
\end{align}
Finally, we specify the requirements on the function $\bg$ and the threshold $\sigma$. We assume in what follows that $\sigma:(0,T)\times \partial \Omega\times \mathbb{R}^2_+  \to \mathbb{R}_+$ is a Carath\'{e}odory mapping such that for almost all $(t,x)\in (0,T)\times \partial \Omega$ and all $(b,\omega)\in \mathbb{R}^2_+$ there holds
\begin{equation}
0\le \sigma(t,x,b,\omega) \le \sigma_{\max}<\infty. \label{maxthres}
\end{equation}
Similarly, we assume that $\bg: (0,T)\times \partial \Omega \times \mathbb{R}^2_+ \times \{\mathbb{R}^3\setminus \b0\}\to \mathbb{R}^3$ is a Carath\'{e}odory mapping. Then, in order to ensure compatibility with the boundary condition \eqref{bc2}, we need to require that, for almost all $(t,x)\in (0,T)\times \partial \Omega$ and all $(b,\omega,\bv)\in  \mathbb{R}^2_+ \times \{\mathbb{R}^3\setminus \b0\}$, there holds
\begin{equation}
\label{compat}
\begin{split}
|\sigma (t,x,\omega,b)|<|\bg(t,x,\omega,b,\bv)|\quad \textrm{ and } \quad \lim_{\bv \to \b0} |\bg(t,x,\omega,b,\bv)| = \sigma (t,x,\omega,b).
\end{split}
\end{equation}
We finish this part by introducing growth, coercivity and further structural assumptions on $\bg$. In order to control the energy of the fluid the natural assumption is that there exists $C_g>0$ such that, for all $\bv\neq \b0$ and all $(\omega,b)$,
\begin{equation}
\bg(t,x,\omega,b,\bv) \cdot \bv \ge -C_g \qquad \textrm{ on } (0,T)\times \partial\Omega. \label{lowg}
\end{equation}
On the other hand, we also need an upper bound on $\bg$ in order to identify the limit in the boundary integrals. Here we assume, roughly speaking,
that the integrability of $\bg \cdot \bv$ (that will be guaranteed by the energy equality and \eqref{lowg}) implies some integrability of $\bg$.  %Therefore, our minimal requirement assumed in the paper is roughly speaking that having the integrability of $\bg \cdot \bv$ (that will be guaranteed by energy equality and \eqref{lowg}) we also have some integrability of $\bg$ that is uniform with respect to $t$ and $x$. To be more precise,
More precisely, we require that there exist $\beta_{g}>1$ and $C>0$ such that,  for all $\bv\neq 0$ and all $(\omega,b)$ and almost all $(t,x)\in (0,T)\times\partial \Omega$,
\begin{equation}
|\bg(t,x,\omega,b,\bv)|^{\beta_{g}}\le C(1+|\bg(t,x,\omega,b,\bv) \cdot \bv| + |\bv(t,x)|^{\frac83}).\label{upg}
\end{equation}
As we claimed before, we also want to cover the possible case when the turbulent energy vanishes on a zero measure set. However, this can cause inability to identify the trace of $\omega$. On the other hand, as it will be clear from the proof, we will always be able to identify the trace of $b\omega$. Therefore, our last assumption on $\bg$ and $\sigma$ is the following. There exist Carath\'{e}odory mappings $\bg^*$ and $\sigma^*$ such that\footnote{We could even assume a more general situation, and adopt the assumption that $\bg$ and $\sigma$ does not depend on $\omega$ whenever $b=0$. However, for the simplicity of the presentation, we do not consider this extension here.}
\begin{equation}\label{lastass}
\bg(t,x,b,\omega,\bv)=\bg^*(t,x,b,b\omega,\bv), \qquad \sigma(t,x,b,\omega)=\sigma^*(t,x,b,b\omega).
\end{equation}

\subsection{Main result}

In order to simplify the presentation of the key result and its proof (but not lose any of the generality of the main theorem) we assume in what follows that all material constants $2\nu_0,\kappa_1,\ldots, \kappa_4$ are equal to one. For the same reason we also introduce $\mu$ to be defined through
\begin{align}
\mu:=\frac{b}{\omega}, \label{dfmunkT1}
\end{align}
and we recall that $E$, the total kinetic energy, is then defined as
\begin{align}
E=\frac{|\bv|^2}{2} + b.\label{dfmunkT2}
\end{align}

\begin{Theorem}
\label{TH1}
Let $\Omega \subset \mathbb{R}^3$ be a $\mathcal{C}^{1,1}$ domain and $T>0$. Assume that the initial data satisfy \eqref{Assv}--\eqref{Asso}, the boundary data satisfy \eqref{bd1}--\eqref{bd4},  $\sigma$ and $\bg$ satisfy \eqref{maxthres}--\eqref{upg} with $\beta_{\bg}>1$. Then, there exists a quintuple $(\bv,b,\omega,p,\bs)$ such that
\begin{align}
\bv &\in L^2(0,T; W^{1,2}_{\bn, \diver}) \cap W^{1,q}(0,T; W^{-1,q}_{\bn}) \textrm{ for all } q\in [1,q_{\min}),\label{k1T}\\
b&\in  \mathcal{E}, \qquad \qquad \qquad \qquad (\mathcal{E} \textrm{ is defined in } \eqref{eee}) \label{k2T}\\
\partial_{t} b &\in \mathcal{M}(0,T; W^{-1,1}_{\Gamma}(\Omega)), \label{k3T}\\
p &\in L^1(0,T; L^1_0(\Omega)), \label{k3P} \\
E &\in W^{1,q}(0,T; W^{-1,q}_{0}(\Omega')) \textrm{ for all } q\in [1,e_{\min}) \textrm{ and all } \overline{\Omega'}\subset \Omega,\label{k2Tx}\\
\omega &\in  L^{\infty}(0,T;L^{\infty}(\Omega)),\label{k4T}\\
\partial_{t} \omega &\in L^q(0,T; W^{-1,q}_{\Gamma}(\Omega))\textrm{ for all } q\in \left[1,\frac{16}{11}\right), \label{k5T}\\
b(\omega-\omega_{\Gamma}) &\in L^q(0,T; W^{1,q}_{\Gamma}(\Omega))\textrm{ for all } q\in \left[1,\frac{16}{11}\right), \label{k5T*}\\
\omega_{\min}e^{-T\omega_{\max}}&\le \omega \le \omega_{\max} \textrm{ a.e. in } Q, \label{minmaxokT}
\end{align}
where
$$
q_{\min}:=\min\left\{\beta_{\bg}, \frac{16}{11}\right\}, \quad \beta_{\min}:=\min \left\{\beta_{\bg},\frac{80}{79}\right\}, \quad \textrm{ and } \quad  E:=\frac{|\bv|^2}{2} + b.
$$
In addition, the pressure $p$ can be decomposed as $p=p_1+p_2+p_3$, where
\begin{align}
p_1 &\in L^{q}(0,T; L^q_0(\Omega))&&\textrm{for all } q\in \left[1, \frac{16}{11}\right), \label{k3P***} \\
p_2 &\in L^{\frac53}(0,T; L^{\frac53}_0(\Omega)), \label{k3P*} \\
p_3 &\in L^{\beta_{\bg}}(0,T; L^{\frac{3\beta_{\bg}}{2}}_0(\Omega)\cap L^{\infty}_{loc}(\Omega))\label{k3P**}
\end{align}
and after denoting
$$
\mu:=\frac{b}{\omega},
$$
the quintuple $(\bv,b,\omega,p,\bs)$ satisfies the following identities:
\begin{align}
&\left. \begin{aligned}
&\int_0^T\langle \partial_{t} \bv, \bw \rangle
-\left(\bv\otimes \bv,
\nabla \bw\right) +(\bs,
\bw)_{\partial \Omega}+\left(\mu\bD(\bv),\bD(\bw)\right) \; dt\\
&\qquad =\int_0^T(p,\diver \bw)\; dt \qquad\textrm{ for all }
\bw \in L^{\infty}(0,T;W^{1,\infty}_{\bn}),
\end{aligned}\right. \label{weak1-1kT}\\
&\left.
\begin{aligned}
&\int_0^T\langle \partial_{t} E, z\rangle -(\bv (E+p), \nabla z) + (\mu \nabla b, \nabla z) + (\mu \bD(\bv) \bv, \nabla z)\; dt\\
&\qquad =-\int_0^T(b \omega, z)\; dt\qquad \textrm{ for all } z\in L^{\infty}(0,T; \mathcal{D}(\Omega)),
\end{aligned}\right. \label{weak3-1kT}\\
&\left.
\begin{aligned}
&\int_0^T\langle \partial_{t} \omega, z\rangle -(\bv\omega, \nabla z) + \left(\frac{\nabla (b\omega)}{\omega} - \nabla b, \nabla z\right)\; dt \\
&\quad = -\int_0^T (\omega^2,z)\; dt\qquad \textrm{ for all } z\in L^{\infty}(0,T;W^{1,\infty}_{\Gamma}(\Omega)),
\end{aligned}\right. \label{weak2-1kT} \\
&\left.
\begin{aligned}
&\bh(\cdot, b, b\omega; \bs, \bv_{\btau}) = \b0 \qquad \textrm{ for almost all } (t,x) \in (0,T)\times \partial \Omega, \\
&\hfill\textrm{ where } \bh \textrm{ is introduced in } \eqref{bc2impl} \textrm{ to describe } \eqref{bc2},
\end{aligned}\right. \label{weak2-ceqs}
\end{align}
with the initial data fulfilling
\begin{equation}
\lim_{t\to 0_+}\|\bv(t)-\bv_0\|_2 + \|\omega(t)-\omega_0\|_2+ \|b(t)-b_0\|_1 =0. \label{oinda333T}
\end{equation}
Moreover, the following inequality holds:
\begin{align}
&\left.\begin{aligned}
&\langle \partial_{t} b, z\rangle +\int_0^T  (\mu \nabla b, \nabla z)-(\bv b, \nabla z)\; dt \ge \int_0^T(-b \omega + \mu|\bD(\bv)|^2, z)\; dt\\
&\qquad \textrm{ for all } z\in \mathcal{C}(0,T;W^{1,\infty}_{\Gamma}(\Omega)) \textrm{ and } z\ge 0 \textrm{ a.e. in } Q.
\end{aligned}\right. \label{entropy1}
\end{align}
In addition, if $\beta_{\bg}>\frac{8}{7}$ then there exists a $\beta_1>1$ such that
\begin{equation}
\partial_t E\in W^{1,\beta_1}(0,T; W^{-1,\beta_1}_0(\Omega))\cap \mathcal{M}(0,T; W^{-1,\beta_1}_{\Gamma}(\Omega)).\label{bettere}
\end{equation}
Moreover, \eqref{weak3-1kT} holds also for all $z\in L^{\infty}(0,T; W^{1,\infty}_{0}(\Omega))$ and the following inequality holds:
\begin{align}
&\left.
\begin{aligned}
&\langle \partial_{t} E, z\rangle +\int_0^T  (\mu \nabla b, \nabla z)-(\bv (E+p), \nabla z)  + (\mu \bD(\bv) \bv, \nabla z) +(b \omega, z)\; dt \\
&\qquad -\int_0^T (\bs, z\bv)_{\partial \Omega}\; dt \le 0\\
&\qquad \textrm{ for all } z\in L^{\infty}(0,T; W^{1,\infty}_{\Gamma}(\Omega)) \textrm{ and } z\ge 0 \textrm{ a.e. in }Q.
\end{aligned}\right. \label{weak3-1kTE}
\end{align}
Further, if $\beta_{\bg}>\frac85$ then there exists a $\beta_2>1$ such that
\begin{equation}
\partial_t E\in W^{1,\beta_2}(0,T; W^{-1,\beta_2}_{\Gamma}(\Omega))\label{bettere2}
\end{equation}
and \eqref{weak3-1kTE} holds with the equality sign.
\end{Theorem}

\subsection{Reduction to One-Equation Model of Turbulence} \label{sec-reduction}

Since the frequency $\omega$ has the dimension $[\textrm{\tt 1/s}]$, the quantity $\frac{\sqrt{b}}{\omega}$ has the dimension of $[\textrm{\tt m}]$ and it can be used as the local measure of the {length scale} of turbulence. Set $\omega = c \frac{\sqrt{b}}{\ell}$. Assuming further that $\ell$ is given, then it follows that the equation for $\omega$ is redundant and $\frac{b}{\omega} = \frac{\ell}{c} \sqrt{b}$ and $b\omega = \frac{c}{\ell} b\sqrt{b}$. Then the Kolmogorov's system \eqref{BM}--\eqref{TKE} reduces to ($k:= \frac{3}{2} b$)
\begin{equation}\label{prandtl}
\begin{split}
\diver \bv &=0, \\
\bv_{,t} + \diver (\bv \otimes \bv) -  \diver \left( \nu(k) \bD(\bv) \right) &= -\nabla p,\\
k_{,t}+\diver (k\bv) - \diver \left(\mu(k) \nabla k\right)& = - \varepsilon(k) + \nu(k) |\bD(\bv)|^2,
    \end{split}
\end{equation}
where
\begin{equation*}
\nu (k) \sim \sqrt{k}, \quad \mu (k) \sim \sqrt{k} \textrm{ and } \varepsilon(k) \sim \sqrt{k} k.
\end{equation*}
This is the model obtained by Prandtl \cite{Prandtl1945}. A general drawback of a one-equation model of turbulence, such as that proposed by Prandtl, is that the length scale of the turbulence has to be known a priori. We refer to \cite{BuLeMa11} for the mathematical theory in the spirit of Theorem \ref{TH1}, for further features related to this system and references regarding the analysis, numerical computations and some applications (further details and more references can be found in a more recent book \cite{ChaLe14}). From the point of view of mathematical analysis of initial and/or boundary-value problems relevant to the Navier-Stokes system with the viscosity depending on other scalar quantity/quantities, we recall several works on analysis of problems related to or motivated by \eqref{prandtl} that were established prior to \cite{ChaLe14}, see \cite{RL94, RL97B, Le97, BCLM1, BCLM2, BCLM3, GaLeLeMuTa03}.

\subsection{History and motivation}\label{history}

In the autumn of 1941 the German army approached Moscow and the Academy of Sciences of U.S.S.R. moved to Kazan. During the workshop of the department of mathematical and physical sciences (January 26-28, 1942), A.~N.~Kolmogorov (1903-1987) introduced his phenomenological two-equation model for the description of the turbulent flows, see the system of governing equations given in \eqref{BM}--\eqref{TKE}. Later on, still in 1942, a brief report of his talk was published in Russian in Izv. Acad. Nauk U.S.S.R., see \cite{Kolmogorov1942}. The report includes the comment that L.~Landau and P.~Kapitsa took part in the discussion. Landau remarked that\footnote{This translation is taken from the Appendix of \cite{Spalding1991}.} \emph{``A.~N. Kolmogorov was the first to provide correct understanding of the local structure of a turbulent flow. As to the equations of turbulent motion, it should be constantly born in mind, in Landau's opinion, that in a turbulent field the presence of rotation in the velocity was confined to a limited region; qualitatively correct equations should lead to just such a distribution of eddies."}.  Kolmogorov's model was unnoticed for almost fifty years. In fact, in 1945 Prandtl \cite{Prandtl1945} proposed a one-equation model of turbulence that can be obtained from Kolmogorov's model assuming that the local length of the turbulence is known a priori, see Subsect. \ref{sec-reduction} above, and as follows from a nice survey article by Spalding, see \cite{Spalding1991}, the one-equation models of turbulence were in place till 1967 when Harlow $\&$ Nakayama proposed the popular k-epsilon model. Brian Spalding, an expert in the field (also an inventor of two-equation models of turbulence), says (\cite{Spalding1991}): \emph{``The fact that for so many years one-equation models continued to be proposed proves \dots that Kolmogorov's two-equation model was indeed a far-from obvious concept"}.  He also links both models referring to Kolmogorov's model: \emph{``The first quantity is actually two-thirds of the turbulence energy (usually given symbol $k$) which appears in currently popular models; and the second, if multiplied by $b$, is proportional to the energy dissipation rate\footnote{In order to
relate Kolmogorov's system to the k-epsilon model we note that ${b}{\omega}$ is proportional to the {energy-dissipation rate} $\epsilon:=2\nu_0 \left\langle \right |\bD(\bv')|^2 \rangle$.} (usually given the symbol  $\varepsilon$ ...) which is the second variable of the model invented much later by Harlow $\&$ Nakayama (1967)."}

Since 1967, several variants of two-equation models of turbulence were proposed, see \cite{Spalding1991} and the books \cite{LS72, MP94, ChaLe14}. Spalding discusses not only different structures of their nonlinearities but he also makes the comment that separates Kolmogorov's model from other two-equation models of turbulence. Spalding \cite{Spalding1991} says: \emph{``The similarity between the Kolmogorov equation and those of the later authors is so great that the one omission is surprising. Kolmogorov made no provision for a source of $\omega$, although he recognized the existence of sink, for which he was bold enough to specify the multiplying constant."} and he continues:
\emph{`` ... The answer is probably that he implicitly presumed that there was a source of $\omega$ in the immediate vicinity of the wall, by fixing the value of $\omega$ there to a finite value by way of a boundary condition."} These specific features and the structure of Kolmogorov's model have served as a  motivation for us to analyze this particular two-equation model of turbulence.

Another motivation for performing the analysis of problems connected with \eqref{BM}--\eqref{TKE} comes from yet another comment of Spalding. He states, see \cite{Spalding1991}: \emph{``~... it is worth pointing out that the question of which of the possible two-equation turbulence models best fits reality has never been seriously investigated."}

\subsection{Highlights}

We conclude this introductory section by summarizing the key features and main difficulties when analyzing \eqref{BM}--\eqref{TKE}.

Several phenomenological two-equation models of turbulence are in place, one of them is a known k-epsilon model, see \cite{MP94}. In this study, we focus on the model proposed by Kolmogorov in 1942 for the following reasons. First, Kolmogorov's model seems to be the first two-equation model of turbulence proposed far in advance of others. Second, as noted by Landau, the formulation, although very brief, is based on Kolmogorov's insight regarding the local structure and properties of turbulent flows. Third, a significant credit to Kolmogorov and his model is given by Spalding who has been an expert in the area for several decades, see \cite{Spalding1991}.

Phenomenological models of turbulence describe flows in terms of averaged quantities (time, spatial or stochastic mean values). It has been conjectured by many scientists, see for example Bardos or Titi \cite{bardos-1998, bardos-titi-2013}, that such flows should be regular. Their conjecture is supported by the analysis of a simplified Smagorinski model of turbulence for which the long-time and large-data \emph{well-posedness} as well as some higher differentiability of the solution are known, see Ladyzhenskaya \cite{Lady67, Lady69} or Pares \cite{pares1992}, while the full regularity (or more precisely even $C^{1,\alpha}$-regularity) is an interesting open question (even when neglecting the inertia or time-derivative of $\bv$). While in Smagorinsky's model the relationship between the Cauchy stress and $\bD(\bv)$ is nonlinear, in Kolmogorov's model the relation between the Cauchy stress tensor and the velocity gradient is \emph{linear}; the generalized viscosity depends however in a specific manner on two scalar quantities $b$ and $\omega$.

The main aim of this study has been to establish long-time and large-data existence theory for Kolmogorov's two-equation model of turbulence in the spirit of Leray \cite{leray34}, Hopf \cite{Hopf51} and Caffarelli, Kohn, Nirenberg \cite{ckn82} (long-time and large-data existence of suitable weak solution). The existence result established here opens the door to the study of regularity properties of such solutions. The scaling of the Navier-Stokes equations plays an important role in the investigation of (partial) regularity associated with the weak solution of the Navier-Stokes equations. Not only does Kolmogorov's system share the same scaling but in fact there is a two-parameter family of scales in which the involved quantities are invariant. More precisely,
if $(\bv,p,\omega,b)$ solves Kolmogorov's system \eqref{BM}-\eqref{TKE}, then, for any $a,b$ and $\theta>0$, the quadruple $(\bv_\theta,p_\theta,\omega_\theta,b_\theta)$, defined through
\begin{align*}
\bv_\theta (t,x)&:= \theta^{a-b} \bv(\theta^{a} t,\theta^b x), \qquad &&p_\theta (t,x):= \theta^{2(a-b)}  p(\theta^{a} t, \theta^b x), \\
\omega_\theta(t,x)&:= \theta^{a} \omega(\theta^{a} t, \theta^b x), \qquad &&b_\theta (t,x):= \theta^{2(a-b)} b(\theta^{a} t, \theta^b x),
\end{align*}
solves Kolmogorov's system as well.

In order to establish the long-time and large-data existence of a suitable weak solution to the initial and boundary-value problem associated with Kolmogorov's PDE system \eqref{BM}-\eqref{TKE} we have to overcome several difficulties which are worth summarizing. First, the measure $\frac{b}{\omega}$ of the effective diffusivity of turbulence and the effective kinematic viscosity could degenerate, which does not allow one to guarantee the integrability of $\nabla \omega$. Using the relation $b\nabla \omega = \nabla (b\omega) - \omega \nabla b$ and the fact that the quantities on the right-hand side are integrable, we found a reformulation of the equation for $\omega$ where we could take the limit. The compactness of $\omega$ is achieved via a variant of the Div-Curl lemma (see \cite{Mu78, Ta78, Ta79}, \cite{feireisl04, FeNo09}). Second, the $L^2$-integrability of $\nabla \bv$ follows from the equation for $b$. Third, the presence of an $L^1$-nonlinearity in \eqref{TKE} is overcome by replacing it by the equation for $b+|\bv|^2/2$, which however requires that the pressure is integrable. The idea applied here goes back to \cite{FeMa06} and \cite{BuFeMa09} where the Navier-Stokes-Fourier system was analyzed. The necessity to have an integrable pressure excludes the no-slip boundary condition from our analysis. We treat, and this is the fourth point worth mentioning, generalized stick-slip boundary conditions. (Note that if the normal traction $\bs$ could be shown to be bounded over $(0,T)\times \partial \Omega$ and if the considered threshold were above the maximal value of $|\bs|$ over $(0,T)\times \partial \Omega$, then the no-slip problem could be successfully analyzed in this way.) Technical difficulties were caused by the fact that we wished to include nonhomogeneous Dirichlet boundary data for $\omega$ and $b$ on part of $\partial \Omega$.

There is an alternative study by Mielke and Naumann, see the announcement of their result in \cite{MiNa15}. Their approach is different in several aspects.
They consider merely the spatially periodic problem, and instead of \eqref{TKE} they only proved the inequality \eqref{TKI}; more precisely they introduce a  nonnegative measure so that the equality holds. They also have a stronger assumption on $b_0$.
In our approach, we investigate flows in bounded domains with the turbulence generated on the boundary. The equivalent formulation of the equation for $b$ proposed here does not require one to introduce a measure into our setting, but requires the integrability of the pressure. We show that an integrable pressure exists even for a very general class of stick-slip boundary conditions. Referring also to \cite[Section 4]{ChaLe14}, we are not aware of any other result concerning long-time and large-data (or well-posedness) existence of (weak) solutions for a two-equation model of turbulence.

\section{Scheme of the proof of Theorem \ref{TH1}}\label{scheme}
The proof of Theorem~\ref{TH1} is constructive and uses a hierarchy of approximations. We introduce them in the following subsections and for each level of approximation we state the result about the existence of solution to the particular approximation. The proofs of these auxiliary lemmas are given in Section~\ref{aupr}. Finally, based on these auxiliary results, we provide the proof of  Theorem~\ref{TH1} in Section~\ref{SMT}.

\subsection{Auxiliary results, inequalities and notations}
We shall first introduce several cut-off functions that will be used when constructing the approximate problems. For any $m\in \mathbb{R}_+$, we define a function $T_m$ as
\begin{equation}
T_m(s):=\left \{\begin{aligned} &s &&\textrm{if } |s|\le m,\\
&m\; \textrm{sgn}\, (s) &&\textrm{if } |s|>m.
\end{aligned}\right.
\end{equation}
We use the symbol $\Theta_m$ to denote the primitive function to $T_m$, i.e.,
\begin{equation}\label{Theta}
\Theta_m(s):=\int_0^s T_m(\tau)\; d\tau.
\end{equation}
Next, we consider a smooth non-increasing function $G$, which is from this point assumed to be fixed, such that
$G(s)=1$ when $s\in [0,1]$ and $G(s)=0$ for $s\ge 2$. Then  for arbitrary $m\in \mathbb{R}_+$ we define
$$
G_m(s):= G\left(\frac{s}{m}\right)
$$
and we denote by $\Gamma_m$ the primitive function to $G_m$, i.e.,
$$
\Gamma_m(s):=\int_0^s G_m(\tau)\; d\tau.
$$
For further purposes we also set $z_{+}:=\max\{0,z\}$ and $z_{-}:=\min\{0,z\}$ for the positive and the negative part of a real number $z$, respectively.

Next, we recall several well-known results from the theory of partial differential equations and function spaces. \emph{Korn's inequality}  states that for any Lipschitz domain $\Omega\subset \mathbb{R}^3$ and any $p\in (1,\infty)$ there exists a constant $C>0$ such that
\begin{equation}
\|\bv\|_{1,p}\le C(\|\bv\|_2 + \|\bD(\bv)\|_p). \label{Korn}
\end{equation}
\emph{The trace theorem} (see \cite[Lemma D.1]{BuGwMaSw12}) states that for a Lipschitz domain $\Omega$, arbitrary $p\in (1,\infty)$ and $\alpha >1/p$, the trace operator is a bounded linear operator from $W^{\alpha,p}(\Omega)$ to $W^{\alpha - 1/p, p}(\partial \Omega)$. In particular, the following estimate holds:
\begin{equation}
\|u\|_{W^{\alpha-\frac{1}{p},p}(\partial \Omega)}\le C\|u\|_{W^{\alpha,p}(\Omega)}.\label{trace}
\end{equation}
We will also require (in order to obtain the optimal estimates of the pressure) $W^{2,p}$ regularity results concerning Poisson's equation with homogeneous Neumann data, i.e., the problem
\begin{equation}\label{laplace}
\begin{aligned}
\triangle u &=g &&\textrm{in } \Omega,\\
\nabla u \cdot \bn &=0 &&\textrm{on } \partial \Omega.
\end{aligned}
\end{equation}
Recalling \cite[Chapter 2]{Gr85}, it is known that for any $\Omega \in \mathcal{C}^{1,1}$, arbitrary $p\in (1,\infty)$ and  $g\in L^p_0(\Omega)$, we can find a unique weak solution of \eqref{laplace} satisfying
\begin{align}
\|u\|_{2,p}&\le C(p,\Omega)\|g\|_p. \label{nabla2}\\
\intertext{In addition, if $g=\diver \bef$ with $\bef\in W^{1,r}_{\bn}$ and $r\in (1,\infty)$, then}
\|u\|_{1,r} & \le C(r,\Omega)\|\bef\|_r. \label{nabla1}
\end{align}

\subsection{$k$-approximation}
The first approximation we introduce here is  an infinite-dimensional $k$-approximation of our problem that will be further approximated by a cascade of finite-dimensional approximations introduced below. Since at the level of the $k$-approximation we want to apply standard monotone operator theory to identify the limit of the last term in \eqref{TKE} (which means that we want to take advantage of the energy equality that comes from the fact that the two formulations of the balance of energy are equivalent for this level of approximation) we use the function $G_k$ to cut the convective term off. Also, in order to avoid difficulties with possibly unbounded turbulent kinetic energy $b$ we cut the viscosity term with the help of the function $T_k$. So the $k$-approximation takes the following form: we want to find $(\bv, p, \omega, b):=(\bv^k, p^k, \omega^k, b^k)$ such that
\begin{align}
\diver \bv &=0, \label{BM-k}\\
\partial_{t} \bv + \diver (G_k(|\bv|^2)\bv \otimes \bv) - \diver \left( T_k\left(\mu \right)\bD(\bv) \right) &= -\nabla p,\label{BLM-k}\\
\partial_{t} \omega + \diver (\omega\bv) - \diver \left (\mu \nabla \omega \right) &= - \omega^2,\label{len-k}\\
\partial_{t} b + \diver (b\bv)-\diver \left(\mu\nabla b\right)&=-b\omega +  T_k\left(\mu\right)|\bD(\bv)|^2,\label{TKE-k}
\end{align}
where $\mu:=\frac{b}{\omega}$. We complete the system with the boundary conditions \eqref{bc1}, \eqref{bc3}--\eqref{bc6} and with the initial conditions \eqref{ID}$_1$ and \eqref{ID}$_3$. In addition, we replace \eqref{bc2} (respectively \eqref{bc2impl}) by the following relation
\begin{equation}
\bg^k (t,x,b,\omega, \bv) + \left(T_k\left(\mu\right)\bD(\bv) \bn\right)_{\btau}=0 \qquad \textrm{on } (0,T)\times \partial \Omega, \label{bc2-k}
\end{equation}
where $\bg^k$ is defined as
\begin{equation}\label{dfgk}
\bg^k(t,x,b,\omega,\bv):= \frac{\bg(t,x,b,\omega,\bv)}{1+k^{-1}|\bg(t,x,b,\omega,\bv)|}\min\{1,k|\bv|\}.
\end{equation}
The reason for this approximation is twofold. First, it is evident that $\bg^k$ is a bounded function (with the bound depending on $k$). Second, we see that due to the presence of $\min\{1,k|\bv|\}$ we can extend continuously $\bg^k$ by zero for $\bv=\b0$.
The last modification is applied to the initial condition for $b$, where we replace $b_0$ by
\begin{equation}
b(0,x)=b_0^k(x):=b_0(x)+\frac{1}{k} \label{ID-k}
\end{equation}
in order to get at this level a proper bound for $b$ from below.

Finally, we neglect the pressure by projecting \eqref{BLM-k} onto the space of divergence-less  test functions and say that for $k\in \mathbb{N}$ fixed, the triple $(\bv, \omega, b)$ solves Problem~$\mathcal{P}_k$ if $(\bv, \omega, b)$ satisfies, in a weak sense, \eqref{BM-k}--\eqref{TKE-k}, completed by the boundary conditions \eqref{bc1}, \eqref{bc3}--\eqref{bc6}, \eqref{bc2-k}--\eqref{dfgk} and initial conditions \eqref{ID}$_1$, \eqref{ID}$_3$ and \eqref{ID-k}. The existence of a weak solution to Problem~$\mathcal{P}_k$ is stated in the following lemma.
\begin{Lem}\label{k-lemma}
Let $\Omega \subset \mathbb{R}^3$ be a Lipschitz domain, $T>0$ be given and let $k\in \mathbb{N}$ fulfilling $k\ge \frac{1}{b_{\min}}$ be arbitrary. Assume that the initial data satisfy \eqref{Assv}--\eqref{Asso}, the boundary data satisfy \eqref{bd1}--\eqref{bd4},  $\bg$ satisfies \eqref{lowg} and \eqref{upg} and $\bg^k$ is defined in \eqref{dfgk}. Then, there exists a triple $(\bv,b,\omega)$ satisfying
\begin{align}
\bv &\in L^2(0,T; W^{1,2}_{\bn, \diver}) \cap W^{1,2}(0,T; W^{-1,2}_{\bn,\diver}),\label{k1}\\
b-b_{\Gamma}&\in  L^q(0,T; W^{1,q}_{\Gamma}(\Omega)) \cap L^{\infty}(0,T; L^1(\Omega)) &&\textrm{for all } q\in \left[1,\frac{5}{4}\right), \label{k2}\\
\partial_{t} b &\in L^1(0,T; W^{-1,1}_{\Gamma}(\Omega)), \label{k3} \\
\omega-\omega_{\Gamma} &\in L^2(0,T; W^{1,2}_{\Gamma}(\Omega))\cap L^{\infty}(0,T;L^{\infty}(\Omega)),\label{k4}\\
\partial_{t} \omega &\in L^q(0,T; W^{-1,q}_{\Gamma}(\Omega))  &&\textrm{for all } q\in \left[1,\frac{16}{11}\right), \label{k5}\\
\omega_{\min}&e^{-T\omega_{\max}}\le \omega \le \omega_{\max} &&\textrm{a.e. in } Q, \label{minmaxok}\\
k^{-1}&e^{-T\omega_{\max}}\le b &&\textrm{a.e. in } Q, \label{minbk}
\end{align}
which solves Problem $\mathcal{P}_k$ in the following sense
\begin{align}
&\left. \begin{aligned}
&\langle \partial_{t} \bv , \bw \rangle
-\left(G_k(|\bv|^2)\bv\otimes \bv,
\nabla \bw\right) + (\bg^k(\cdot,b,\omega,\bv),
\bw)_{\partial \Omega} \\
&\quad +\left(T_k(\mu)\bD(\bv),\bD(\bw)\right) =0
\qquad \textrm{ for all }
\bw \in W^{1,2}_{\bn,\diver} \textrm{ and a.a. } t\in (0,T),
\end{aligned}\right. \label{weak1-1k}\\
&\left.
\begin{aligned}
&\langle \partial_{t} b, z\rangle -(b \bv, \nabla z) + (\mu \nabla b, \nabla z) = (-b \omega + T_k(\mu)|\bD(\bv)|^2, z)\\
&\qquad \textrm{ for all } z\in W^{1,\infty}_{\Gamma}(\Omega) \textrm{ and a.a. } t\in (0,T),
\end{aligned}\right. \label{weak3-1k}\\
&\left.
\begin{aligned}
&\langle \partial_{t} \omega, z\rangle -(\omega \bv, \nabla z) + (\mu \nabla \omega, \nabla z) = -(\omega^2,z)\\ &\qquad \textrm{ for all } z\in W^{1,\infty}_{\Gamma}(\Omega) \textrm{ and a.a. } t\in (0,T),
\end{aligned}\right. \label{weak2-1k}
\end{align}
where $\mu$ is given as
\begin{align}
\mu&:=\frac{b}{\omega}.\label{dfmunk}
\end{align}
The initial data are attained strongly in the corresponding spaces, i.e.,
\begin{equation}
\lim_{t\to 0_+}\|\bv(t)-\bv_0\|_2 + \|\omega(t)-\omega_0\|_2+ \|b(t)-b^{k}_0\|_1 =0. \label{oinda333}
\end{equation}
Moreover, for all $\lambda \in (0,1]$ the following uniform ($k$-independent) estimate holds
\begin{equation}
\begin{split}
&\sup_{t\in (0,T)} \left(\|b(t)\|_1 + \|\ln b(t)\,\|_1+\|\bv(t)\|^2_2\right) + \int_Q (1+b^{-1})T_k(\mu)|\bD(\bv)|^2 \; dx \; dt\\
&\quad + \int_0^T \left( \int_{\Omega} \frac{\mu}{b^{1+\lambda}}|\nabla b|^2+\mu |\nabla \omega|^2 +|\mu|^{\frac{8}{3}-\lambda}\; dx + \int_{\partial \Omega} |\bg^k\cdot \bv|\; dS \; \right) dt\\
&\quad \le C(\lambda^{-1},\bv_0,b_0,\omega_0, \omega_{\min},\omega_{\max},b_{\min},b_{\max}).
\end{split}
\label{finalap}
\end{equation}
\end{Lem}

\subsection{$(n,k)$-approximation}
In order to prove Lemma~\ref{k-lemma} we use a Galerkin approximation for the velocity to replace \eqref{weak1-1k}. Moreover, since we want to use the standard $L^2$-theory for $b$ and $\omega$, we add to $\mu$ the coefficient $1/n$ and replace $\mu$ by $T_n(\mu)$ in all diffusion terms in  \eqref{len-k}--\eqref{TKE-k} and we also replace the term on the right-hand side of \eqref{TKE-k} by its proper truncation (see below). Moreover, we mollify the initial condition $b_0$ in the following way. We find a sequence $\{b_0^n\}_{n=1}^\infty$ of smooth nonnegative functions such that
\begin{equation}
b_0^n \to b_0 \qquad \textrm{strongly in } L^1(\Omega) \label{bcon2}
\end{equation}
and consider now  the initial condition
\begin{equation}
b_0^{n,k}:= b_0^n + \frac{1}{k}.\label{inn}
\end{equation}
In addition, we also mollify the boundary data $b_{\Gamma}$ as follows. We find a sequence $\{b^n_{\Gamma}\}_{n=1}^\infty$ of smooth functions satisfying \eqref{bd4} such that
\begin{align}\label{bboundar}
b_{\Gamma}^n &\to b_{\Gamma} &&\textrm{strongly in } L^2(0,T; W^{1,2}(\Omega))\cap W^{1,1}(0,T; L^1(\Omega)).
\end{align}

To summarize,  let  $\{\bw_i\}_{i=1}^\infty$ be a basis of $W^{1,2}_{\bn,\diver}$ that is orthogonal in $L^2(\Omega)^3$ (such basis can be easily constructed e.g., by taking the eigenfunctions of the Stokes operator subjected to Neumann boundary conditions) and denote by $V^n$ the linear span of $\{\bw_i\}_{i=1}^n$. We further project the initial data for $\bv$ to the space $V^n$ and denote
\begin{equation}
\bv_0^n:=\sum_{i=1}^{n}c_i^0\bw_i \quad \textrm{ where } \quad c_i^0:=(\bv_0,\bw_i).\label{v0n}
\end{equation}
Note that it follows from \eqref{v0n} that
\begin{align}
\bv_0^n& \to \bv_0 &&\textrm{strongly in } L^2_{\bn, \diver}. \label{conv0}
\end{align}
We shall refer to the problem described above as Problem $\mathcal{P}_{n,k}$ and we  state the existence result for this problem in the following lemma.
\begin{Lem}\label{n-ap}
Let $\Omega \subset \mathbb{R}^3$ be a Lipschitz domain, $T>0$ be given and let $k,n\in \mathbb{N}$ fulfilling $k\ge \frac{1}{b_{\min}}$ be arbitrary. Assume that the initial data satisfy \eqref{Assv}--\eqref{Asso}, the boundary data satisfy \eqref{bd1}--\eqref{bd4},  $\bg$ satisfies \eqref{lowg} and \eqref{upg} and $\bg^k$ is defined in \eqref{dfgk}. Then, there exists a triple $(\bc,b,\omega):=(\bc^n, b^n, \omega^n)$ satisfying
\begin{align}
\bc &\in W^{1,\infty}(0,T)^n,\label{cm1df}\\
b-b_{\Gamma}^n&\in  L^2(0,T; W^{1,2}_{\Gamma}(\Omega)) \cap L^{\infty}(0,T; L^2(\Omega)), \label{cm4df}\\
\partial_{t} b &\in L^2(0,T; W^{-1,2}_{\Gamma}(\Omega)), \label{cm41df}\\
\omega-\omega_{\Gamma} &\in L^2(0,T; W^{1,2}_{\Gamma}(\Omega))\cap L^{\infty}(0,T;L^{\infty}(\Omega)),\label{cm5df}\\
\partial_{t} \omega &\in L^2(0,T; W^{-1,2}_{\Gamma}(\Omega)), \label{cm5.1df}\\
\omega_{\min}e^{-T\omega_{\max}}&\le \omega \le \omega_{\max} \qquad \,\,\textrm{ a.e. in } Q, \label{minmaxon}\\
k^{-1}e^{-T\omega_{\max}}&\le b \qquad \qquad \qquad \textrm{ a.e. in } Q, \label{minbn}
\end{align}
which solves Problem $\mathcal{P}_{k,n}$ in the following sense:
\begin{align}
&\left. \begin{aligned}
&(\partial_{t} \bv, \bw_i )
-\left(G_k(|\bv|^2)\bv\otimes \bv,
\nabla \bw_i\right) + (\bg^k(\cdot,b,\omega,\bv),
\bw_i)_{\partial \Omega} \\
&\qquad +\left(T_k(\mu^n)\bD(\bv),\bD(\bw_i)\right) =0
\qquad \textrm{ for all }
i=1,\ldots, n,
\end{aligned}\right. \label{weak1-1m}\\
&\left.
\begin{aligned}
&\langle \partial_{t} b, z\rangle -(b\bv, \nabla z) + (T_n(\mu^n) \nabla b, \nabla z) = \left(-b \omega + \frac{T_k(\mu^n)|\bD(\bv)|^2}{1+n^{-1}|\bD(\bv)|^2}, z\right)\\
&\qquad \textrm{ for all } z\in W^{1,2}_{\Gamma}(\Omega) \textrm{ and a.a. } t\in (0,T),
\end{aligned}\right. \label{weak3-1m}\\
&\left.
\begin{aligned}
&\langle \partial_{t} \omega, z\rangle -(\omega\bv, \nabla z) + (T_n(\mu^n) \nabla \omega, \nabla z) = -(\omega^2,z)\\
&\qquad \textrm{ for all } z\in W^{1,2}_{\Gamma}(\Omega) \textrm{ and a.a. } t\in (0,T),
\end{aligned}\right. \label{weak2-1m}
\end{align}
with $\bv$ and $\mu^n$ defined as
\begin{align}
\bv(t,x)&:=\sum_{i=1}^n c_i(t)\bw_i(x),\label{vn}\\
\mu^n&:=\frac{b}{\omega}+ \frac{1}{n}.\label{dfmunn}
\end{align}
The initial data are attained in the following sense
\begin{equation}
\lim_{t\to 0_+}\|\bv(t)-\bv^n_0\|_2 + \|\omega(t)-\omega_0\|_2+ \|b(t)-b^{n,k}_0\|_2 =0. \label{oinda2}
\end{equation}

\end{Lem}

\subsection{$(m,n,k)$-Galerkin approximation}
The next approximation we introduce here consists in projecting \eqref{TKE-k} onto a finite dimensional space. Since at the level of Galerkin approximations we do not control the sign of $b$, we also replace  $b$ by its positive part $b_+:= \max\{0,b\}$ in some terms. In addition,  we redefine $\mu$ so that it can not blow-up for singular $\omega$. To be more specific,  let $\{z_i\}_{i=1}^{\infty}$ be a basis of $W^{1,2}_{\Gamma}(\Omega)$ that is orthogonal in $L^2(\Omega)$ and denote by $Z^m$ the linear span of $\{z_i\}_{i=1}^m$. The $m$-dimensional projection of the initial condition $b_0^{n,k}$ is then given as
\begin{equation}
b_0^{m,n,k}:= \sum_{i=1}^m d^0_i z_i+ b^n_{\Gamma}(0),\textrm{ where } d_i^0:=(b_0^{n,k}-b_{\Gamma}^n(0),z_i).\label{b0mnk}
\end{equation}
Note that \eqref{b0mnk} implies that
\begin{align}
b_0^{n,m,k} &\to b_0^{n,k} &&\textrm{strongly in } L^2(\Omega).\label{conb0}
\end{align}
Moreover, to avoid an additional approximation, we mollify the boundary condition for $\omega$, i.e., we find a sequence of smooth functions $\{\omega_{\Gamma}^m\}_{m=1}^{\infty}$ satisfying \eqref{bd3} such that
\begin{align}\label{wboundar}
\omega_{\Gamma}^m &\to \omega_{\Gamma} &&\textrm{strongly in } L^2(I; W^{1,2}(\Omega))\cap W^{1,1}(I; L^1(\Omega)).
\end{align}
The following lemma states the existence of the solution to the problem described in this subsection that we denote Problem $\mathcal{P}_{m,n,k}$.
\begin{Lem}\label{m-ap}
Let $\Omega \subset \mathbb{R}^3$ be a Lipschitz domain, let $T>0$ be given and let $k,m,n\in \mathbb{N}$ fulfilling $k\ge \frac{1}{b_{\min}}$ and $m\ge \omega_{\max}$ be arbitrary. Assume that the initial data satisfy \eqref{Assv}--\eqref{Asso}, the boundary data satisfy \eqref{bd1}--\eqref{bd4},  $\bg$ satisfies \eqref{lowg} and \eqref{upg} and $\bg^k$ is defined in \eqref{dfgk}. Then, there exists a triple $(\bc,\bd,\omega)$ satisfying
\begin{align}
\bc& \in W^{1,\infty}(0,T)^n, \label{cm}\\
\bd&\in W^{1,\infty}(0,T)^m, \label{dm}\\
\omega -\omega_{\Gamma}^m &\in L^2(I; W^{1,2}_{\Gamma}(\Omega)), \label{omm}\\
\partial_{t} \omega&\in L^2(I;W^{-1,2}_{\Gamma}(\Omega)), \label{tdomm}\\
\omega_{\min}e^{-T\omega_{\max}}&\le \omega \le \omega_{\max} \textrm{ a.e. in } Q, \label{minmaxom}
\end{align}
which solves Problem $\mathcal{P}_{m,n,k}$ in the following sense
\begin{align}
&\left. \begin{aligned}
&(\partial_{t} \bv, \bw_i )
-\left(G_k(|\bv|^2)\bv\otimes \bv,
\nabla \bw_i\right) + (\bg^k(\cdot,b,\omega,\bv),
\bw_i)_{\partial \Omega} \\
&\qquad +\left(T_k(\mu^{n,m})\bD(\bv),\bD(\bw_i)\right) =0
\qquad \textrm{ for all }
i=1,\ldots, n,
\end{aligned}\right. \label{weak1-1l}\\
&\left.
\begin{aligned}
&(\partial_{t} b, z_i) -(b \bv, \nabla z_i) + (T_n(\mu^{n,m}) \nabla b, \nabla z_i) \\
&\quad = \left(-b_+ \omega + \frac{T_k(\mu^{n,m})|\bD(\bv)|^2}{1+n^{-1}|\bD(\bv)|^2}, z_i\right)
\qquad \textrm{ for all } i=1,\ldots, m,
\end{aligned}\right. \label{weak3-1l}\\
&\left.
\begin{aligned}
&\langle \partial_{t} \omega, z\rangle -(\omega \bv, \nabla z) + (T_n(\mu^{n,m}) \nabla \omega, \nabla z) = -(T_m(\omega) \omega_+,z)\\ &\qquad \textrm{ for all } z\in W^{1,2}_{\Gamma}(\Omega) \textrm{ and a.a. } t\in (0,T),
\end{aligned}\right. \label{weak2-1l}
\end{align}
with $\mu^{n,m}$, $\bv$ and $b$  given as
\begin{align}
\bv(t,x)&:=\sum_{i=1}^n c_i(t)\bw_i(x),\label{vnm}\\
b(t,x)&:=\sum_{i=1}^m d(t)z_i(x) + b_{\Gamma}^n(t,x),\label{bnm}\\
\mu^{n,m}&:=\frac{b_+}{\omega + \frac{1}{m}}+ \frac{1}{n}.\label{dfmunm}
\end{align}
The initial data are attained in the following sense
\begin{equation}
\lim_{t\to 0_+} \|b(t)-b_0^{n,m,k}\|_2 + \|\bv(t)-\bv_0^n\|_2 + \|\omega(t) - \omega_0\|_2 =0, \label{atainm}
\end{equation}
where $\bv_0^n$ is given in \eqref{v0n} and $b_0^{n,m,k}$ is defined in \eqref{b0mnk}.
\end{Lem}

\section{Proof of auxiliary existence results}\label{aupr}
In this section we shall prove all auxiliary assertions stated in Section~\ref{scheme}, i.e., Lemmas~\ref{k-lemma}--\ref{m-ap}.
\subsection{Proof of Lemma~\ref{m-ap}}\label{sub3.1}
Assume that $m,n,k\in \mathbb{N}_+$ are fixed, $k\ge \frac{1}{b_{\min}}$ and $m\ge \omega_{\max}$, and that all assumptions of Lemma~\ref{m-ap} are satisfied. Recall that $\{z_i\}_{i=1}^{\infty}$ and $\{\bw_{i}\}_{i=1}^{\infty}$ denotes the basis of $W^{1,2}_{\Gamma}(\Omega)$ and $W^{1,2}_{\bn,\diver}$ orthogonal in $L^2(\Omega)$ and $L^2_{\bn,\diver}$, respectively. To prove Lemma~\ref{m-ap}, we consider the Galerkin approximation of the last equation which has not been approximated yet, i.e.,  for given arbitrary $\ell \in \mathbb{N}_+$ we look for $(\bv^{\ell},\omega^{\ell},b^{\ell})$ given as
\begin{align}
\bv^{\ell}(t,x)&:=\sum_{i=1}^n c^{\ell}_i(t)\bw_i(x),\label{vnml}\\
b^{\ell}(t,x)&:=\sum_{i=1}^m d^{\ell}(t)z_i(x) + b_{\Gamma}^n(t,x),\label{bnml}\\
\omega^{\ell}(t,x)&:=\sum_{i=1}^{\ell} e_i^{\ell}(t) z_i(x)+ \omega_{\Gamma}^m(t,x),\label{onml}
\end{align}
and we require that the coefficients $\bc^{\ell}=(c_1^{\ell},\ldots, c_n^{\ell})$, $\bd^{\ell}=(d^{\ell}_1,\ldots, d^{\ell}_m)$ and $\be^{\ell}=(e_1^{\ell}, \ldots, e_{\ell}^{\ell})$ solve  the following system of ordinary differential equations on $(0,T)$:
\begin{align}
&\left. \begin{aligned}
&(\partial_{t} \bv^{\ell}, \bw_i )
-\left(G_k(|\bv^{\ell}|^2)\bv^{\ell}\otimes \bv^{\ell},
\nabla \bw_i\right) + (\bg^k(\cdot, b^{\ell}, \omega^{\ell}, \bv^{\ell}),
\bw_i)_{\partial \Omega} \\
&\qquad +\left(T_k(\mu^{\ell})\bD(\bv^{\ell}),\bD(\bw_i)\right) =0\qquad \textrm{ for all }
i=1,\ldots, n,
\end{aligned}\right. \label{weak1-1}\\
&\left.
\begin{aligned}
&(\partial_{t} b^{\ell}, z_i) -(b^{\ell} \bv^{\ell}, \nabla z_i) + (T_n(\mu^{\ell}) \nabla b^{\ell}, \nabla z_i)\\
 &\qquad = \left(-b^{\ell}_+ T_m(\omega^{\ell}_+) + \frac{T_k(\mu^{\ell})|\bD(\bv^{\ell})|^2}{1+n^{-1}|\bD(\bv^{\ell})|^2}, z_i\right) \qquad \textrm{ for all } i=1,\ldots, m,
\end{aligned}\right. \label{weak3-1}\\
&\left.
\begin{aligned}
&(\partial_{t} \omega^{\ell}, z_i) -(\omega^{\ell} \bv^{\ell}, \nabla z_i) + (T_n(\mu^{\ell}) \nabla \omega^{\ell}, \nabla z_i) = -(T_m(\omega^{\ell})\omega^{\ell}_+,z_i) \\
&\qquad \textrm{ for all } i=1,\ldots, \ell,
\end{aligned}\right. \label{weak2-1}
\end{align}
where
\begin{equation}
\mu^{\ell}:=\frac{b^{\ell}_+}{\omega^{\ell}_+ + \frac{1}{m}}+ \frac{1}{n}.\label{dfmunml}
\end{equation}
We consider the initial conditions for $(\bc^{\ell},\bd^{\ell},\be^{\ell})$ described by the following relations
\begin{equation}
\begin{split}
\bv^{\ell}(0)&= \bv_0^n, \quad
b^{\ell}(0)= b_0^{m,n,k}, \quad
\omega^{\ell}(0)= \omega_0^{\ell},
\end{split}
\label{indata1}
\end{equation}
where $\bv_0^n$ and $b_0^{m,n,k}$ are defined in \eqref{v0n} and \eqref{b0mnk} respectively, and
\begin{equation}\label{omega0}
\omega_0^{\ell}:=\sum_{i=1}^{\ell} e_i^0 z_i +\omega_{\Gamma}^m(0) \quad \textrm{ with } \quad e_i^0:= (\omega_0-\omega^m_{\Gamma}(0), z_i).
\end{equation}
Note that it directly follows from this definition that
\begin{align}
\omega_0^{\ell} \to \omega_0 &&\textrm{strongly in } L^2(\Omega). \label{cono0}
\end{align}
The existence of a solution to \eqref{bnml}--\eqref{indata1} on a short time interval follows from Carath\'{e}odory's theorem. Moreover, using the a~priori estimates (independent of $\ell$ and $t\in (0,T)$) established below, we can extend the solution onto the whole time interval $(0,T)$. Our goal is to let $\ell \to \infty$ to obtain the statement of Lemma~\ref{m-ap}, which we will do in the following subsections.

\subsubsection{Uniform $\ell$-independent estimates} Here, and also in what follows, we use a generic  constant $C$ to indicate independence of a quantity on $(k,n,m,\ell)$. If some estimates depend on some parameters this will be clearly indicated in the text.

First, multiplying the $i$th equation in \eqref{weak1-1} by $c_i^{\ell}$ and summing the result over $i=1,\ldots, n$ we get the identity
\begin{equation}
\begin{split}
&\frac12 \frac{d}{dt}\|\bv^{\ell}\|_2^2 -
\frac12(G_k(|\bv^{\ell}|^2)\bv^{\ell}, \nabla |\bv^{\ell}|^2) +\left(\bg^k(\cdot,\omega^{\ell},b^{\ell},\bv^{\ell}),\bv^{\ell}\right)_{\partial \Omega} \\
&\quad + \int_{\Omega}T_k(\mu^{\ell})|\bD(\bv^{\ell})|^2 \; dx =0.
\end{split}\label{apeq}
\end{equation}
Next, using the facts that $\bv^{\ell}\cdot \bn=0$ on $(0,T)\times \partial\Omega$ and $\diver \bv^{\ell}=0$ in $Q$ we deduce that
$$
\frac12(G_k(|\bv^{\ell}|^2)\bv^{\ell}, \nabla |\bv^{\ell}|^2)=\frac12(\bv^{\ell},
\nabla \Gamma_k(|\bv^{\ell}|^2))=-\frac12 (\diver
\bv^{\ell},\Gamma_k(|\bv^{\ell}|^2))=0.
$$
Thus, using \eqref{apeq},  the nonnegativity of $\mu^{\ell}$, the assumption \eqref{lowg} combined with \eqref{dfgk} and \eqref{conv0} we get
\begin{equation}
\sup_{t\in(0,T)}\|\bv^{\ell}(t)\|_2^2 + \int_Q T_k(\mu^{\ell})|\bD(\bv^{\ell})|^2\; dx \; dt \le C. \label{apest1}
\end{equation}
Here, $C$ is greater than $2 C_{\bg} T|\partial\Omega| + \|\bv_0\|_2^2$.
Consequently, using the orthonormality of the basis $\{\bw_i\}$ in $L^2(\Omega)$ we deduce from \eqref{apest1} that
\begin{equation}
\sup_{t\in (0,T)}|\bc^{\ell}(t)|\le C.\label{apestsim}
\end{equation}
Then, using \eqref{weak1-1},  \eqref{apestsim}, the  fact that $|\bg^k|\le k$ and the above estimate \eqref{apestsim}, we can easily obtain
\begin{equation}
\sup_{t\in (0,T)}|\partial_{t} \bc^{\ell}(t)|\le C(n,k,m).\label{apestsim2}
\end{equation}

Next, multiplying the $i$th equation in \eqref{weak3-1} by $d_i^{\ell}$ and summing the result over $i=1,\ldots, m$ (which means that $b^{\ell}-b_{\Gamma}^n$ appears as a ``test function" in \eqref{weak3-1}) we get the identity
\begin{equation}
\begin{split}
(\partial_{t} b^{\ell}, b^{\ell}-b_{\Gamma}^n) - (b^{\ell} \bv^{\ell} , \nabla (b^{\ell}-b^n_{\Gamma})) + (T_n(\mu^{\ell})\nabla b^{\ell}, \nabla (b^{\ell} - b_{\Gamma}^n))\\
=-\left(b_+^{\ell} T_m(\omega^{\ell}_{+}) + \frac{T_k(\mu^{\ell})|\bD(\bv^{\ell})|^2}{1+n^{-1}|\bD(\bv^{\ell})|^2}, b^{\ell}-b_{\Gamma}^n\right).
\end{split}
\end{equation}
Hence, using the smoothness of $b_{\Gamma}^n$, the fact the $\diver \bv^{\ell}=0$ and \eqref{apest1}, we deduce with the help of Gronwall's lemma and Young's inequality that
\begin{equation}
\sup_{t\in (0,T)}\|b^{\ell}\|_2^2 \le C(n,k,m). \label{wel}
\end{equation}
Similarly as before, using the equivalence of norms on finite-dimensional spaces and the identity \eqref{weak3-1} we find that
\begin{equation}
\sup_{t\in (0,T)}(|\bd^{\ell}(t)| + |\partial_{t} \bd^{\ell}(t)|) \le C(n,k,m). \label{wel2}
\end{equation}

Finally, we derive uniform estimates for $\omega^{\ell}$. Multiplying the $i$th equation in \eqref{weak2-1} by $e_i^{\ell}$ and summing over $i=1,\ldots ,{\ell}$ we get the identity
\begin{equation}
\begin{split}
(\partial_{t} \omega^{\ell}, \omega^{\ell}-\omega^m_{\Gamma}) - (\omega^{\ell} \bv^{\ell}, \nabla (\omega^{\ell}-\omega^m_{\Gamma})) + (T_n(\mu^{\ell})\nabla \omega^{\ell}, \nabla (\omega^{\ell} - \omega^m_{\Gamma})\\
=-(T_m(\omega^{\ell})\omega_+^{\ell}, \omega^{\ell}-\omega^m_{\Gamma}). \label{ID1o}
\end{split}
\end{equation}
Consequently, adding and subtracting terms with $\omega_{\Gamma}^m$ to the corresponding integrals and  using the divergence-free constraint on $\bv^{\ell}$, we deduce from \eqref{ID1o} with the help of integration by parts that
\begin{equation}
\begin{split}
&\frac12 \frac{d}{dt} \| \omega^{\ell}-\omega_{\Gamma}^m\|_2^2  + \int_{\Omega} T_n(\mu^{\ell})|\nabla (\omega^{\ell}-\omega_{\Gamma}^m)|^2\; dx+(T_m(\omega^{\ell})\omega_+^{\ell}, \omega^{\ell}-\omega^m_{\Gamma})\\
& = -(T_n(\mu^{\ell})\nabla \omega^m_{\Gamma}, \nabla (\omega^{\ell} - \omega_{\Gamma}^m))-(\partial_{t} \omega^m_{\Gamma},\omega^{\ell}-\omega^m_{\Gamma}) -((\omega^{\ell}-\omega^m_{\Gamma})\bv^{\ell}, \nabla  \omega^m_{\Gamma}). \label{ID1o2}
\end{split}
\end{equation}
Thus, using the definition of $\mu^{\ell}$ (see \eqref{dfmunml}),   the properties of the function $T_n$, the bound \eqref{apest1}, the smoothness of $\omega_{\Gamma}^m$ and Young's inequality, we find that
\begin{equation}
\begin{split}
\frac{d}{dt} \| \omega^{\ell}-\omega_{\Gamma}^m\|_2^2  + \frac{1}{n}\|\nabla (\omega^{\ell}-\omega_{\Gamma}^m)\|_2^2 \le C(n,m,k) (\|\omega^{\ell}-\omega_{\Gamma}^m)\|_2^2 + 1). \label{ID1o3}
\end{split}
\end{equation}
Hence, we see that by using  Gronwall's lemma we get the uniform ($\ell$-independent estimate)
\begin{equation}
\begin{split}
\sup_{t\in (0,T)} \| \omega^{\ell}(t)-\omega_{\Gamma}^m(t)\|_2^2  + \int_0^T \|\nabla (\omega^{\ell}-\omega_{\Gamma}^m)\|_2^2 \; dt \le C(n,m,k). \label{apestom1}
\end{split}
\end{equation}
Having \eqref{apestom1}, it is then standard to deduce from \eqref{weak2-1} that
\begin{equation}
\int_0^T \|\partial_{t} \omega^{\ell}\|^2_{W^{-1,2}_{\Gamma}}\; dt \le C(n,m,k). \label{otim}
\end{equation}

\subsubsection{Taking the limit $\ell \to \infty$}
Using \eqref{apestsim}--\eqref{apestsim2} and \eqref{wel2} we can find a subsequence (that we do not relabel) such that
\begin{align}
\bc^{\ell} &\rightharpoonup^* \bc &&\textrm{weakly$^*$ in } W^{1,\infty}(0,T)^n,\label{cl1}\\
\bd^{\ell} &\rightharpoonup^* \bd &&\textrm{weakly$^*$ in } W^{1,\infty}(0,T)^m. \label{cl2}\\
\intertext{Consequently, using the Arsela-Ascoli theorem and the definition of $\bv^{\ell}$ and $b^{\ell}$, we see that}
\bc^{\ell} &\to \bc &&\textrm{strongly in } \mathcal{C}(0,T)^n,\label{cl1s}\\
\bd^{\ell} &\to \bd &&\textrm{strongly in } \mathcal{C}(0,T)^m, \label{cl2s}\\
\bv^{\ell}&\to \bv &&\textrm{strongly in } \mathcal{C}(0,T; W^{1,2}_{\bn, \diver}), \label{cl3}\\
b^{\ell}-b_{\Gamma}^n&\to b-b_{\Gamma}^n &&\textrm{strongly in } \mathcal{C}(0,T; W^{1,2}_{\Gamma}). \label{cl4}
\end{align}
Moreover, using \eqref{apestom1}--\eqref{otim} and the Aubin-Lions lemma, we can find a subsequence (that is again not relabelled) such that
\begin{align}
\omega^{\ell}-\omega_{\Gamma}^m &\rightharpoonup \omega-\omega_{\Gamma}^m &&\textrm{weakly in } L^2(0,T; W^{1,2}_{\Gamma}(\Omega)),\label{cl5}\\
\partial_{t} \omega^{\ell} &\rightharpoonup \partial_{t} \omega &&\textrm{weakly in } L^2(0,T; W^{-1,2}_{\Gamma}(\Omega)), \label{cl5.1}\\
\omega^{\ell} &\to \omega &&\textrm{strongly in } L^2(0,T; L^2(\Omega)). \label{cl6}
\end{align}
Having the convergence results \eqref{cl1}--\eqref{cl6}, it is easy to identify the limit of \eqref{vnml}--\eqref{bnml} and get \eqref{vnm}--\eqref{bnm}. In addition, it is also quite standard to take the limit in   \eqref{weak1-1}--\eqref{weak2-1} and in \eqref{dfmunml} and obtain \eqref{weak1-1l}--\eqref{weak2-1l} and \eqref{dfmunm}, provided that we show that the limit $\omega$ satisfies \eqref{minmaxom}, which we shall show next. (Note that this is the reason why we assumed $m\ge \omega_{\max}$.) The attainment of the initial data \eqref{atainm} can be proven by standard arguments.
\subsubsection{Minimum and maximum principle for $\omega$}
It remains to show \eqref{minmaxom}. To do so, we first identify the limit of \eqref{weak2-1} (without assuming the validity of \eqref{minmaxom}) and get
\begin{align}
&\left.
\begin{aligned}
&\langle \partial_{t} \omega, z\rangle -(\omega\bv, \nabla z) + (T_n(\tilde{\mu}^{n,m}) \nabla \omega, \nabla z) = -(\omega^2,z)\\
&\qquad \textrm{ for all } z\in W^{1,2}_{\Gamma}(\Omega) \textrm{ and a.a. } t\in (0,T),
\end{aligned}\right. \label{weak2-1l*}
\end{align}
where $\tilde{\mu}^{n,m}$ is given as
\begin{equation}
\tilde{\mu}^{n,m}:=\frac{b_+}{\omega_+ + \frac{1}{m}}+ \frac{1}{n}.\label{dfmunm*}
\end{equation}
Since $\omega^m_{\Gamma}\ge \omega_{\min}$ we see that $\omega_- \in L^2(I; W^{1,2}_{\Gamma}(\Omega))$ is a possible test function in \eqref{weak2-1l*} (note that $\omega_{-}\le 0$ and $\omega_{-}\in W^{1,2}_{\Gamma}(\Omega)$ for almost all $t\in (0,T)$. Thus, considering $z=\omega_-$ we observe that the term on the right-hand side of \eqref{weak2-1l*} is identically zero. In addition, using integration by parts we also find (due to the fact that $\diver \bv =0$ and $\bv \cdot \bn =0$ on $\partial \Omega$) that the second term in \eqref{weak2-1l*} is zero. Finally since $\tilde\mu^{n,m}\ge 0$ we observe that
$$
\frac{d}{dt} \|\omega_-\|_2^2 \le 0.
$$
Consequently, using \eqref{Asso}, we conclude that
\begin{equation}
\omega \ge 0 \qquad \textrm{ a.e. in }Q
 \label{oge0}
\end{equation}
and we see that we can replace $\omega_+$ by $\omega$ in \eqref{weak2-1l*} and \eqref{dfmunm*} (hence $\tilde{\mu}^{n,m} = \mu^{n,m}$). Similarly, setting $z:=(\omega - \omega_{\max})_+$ in \eqref{weak2-1l*} (which is again an admissible test function since $\omega_{\Gamma}^m \le \omega_{\max}$) we find by using the same procedure as above (note that the term on the right-hand side of  \eqref{weak2-1l*} is non-positive) and by using \eqref{Asso} that
\begin{equation}
\omega \le \omega_{\max} \qquad \textrm{ a.e. in }Q. \label{ole}
\end{equation}
Since we assume that $m\ge \omega_{\max}$ we can replace $T_m(\omega)$ by $\omega$ in \eqref{weak2-1l*} and conclude that \eqref{weak2-1l*} leads to \eqref{weak2-1l}.
Finally, we set $z:=e^{\omega_{\max}t}(\omega e^{\omega_{\max} t}-\omega_{\min})_{-}$
in \eqref{weak2-1l} (which is again admissible). Note that the convective term again vanishes and the third term on the left-hand side generates a nonnegative term. Thus, we obtain the following inequality
\begin{equation}
\langle \partial_{t} \omega, e^{\omega_{\max}t}(\omega e^{\omega_{\max} t}-\omega_{\min})_{-}\rangle\le -(\omega^2, e^{\omega_{\max}t}(\omega e^{\omega_{\max} t}-\omega_{\min})_{-}) \label{minp2}
\end{equation}
for almost all $t\in (0,T)$. Since,
\begin{equation*}
\begin{split}
\langle \partial_{t} \omega, e^{\omega_{\max}t}(\omega e^{\omega_{\max} t}-\omega_{\min})_{-}\rangle &= \frac12 \frac{d}{dt}\|(\omega e^{\omega_{\max} t}-\omega_{\min})_{-}\|_2^2\\
 &\quad - \omega_{\max}(\omega e^{\omega_{\max}t},(\omega e^{\omega_{\max} t}-\omega_{\min})_{-})
 \end{split}
\end{equation*}
we get from \eqref{minp2} that
\begin{equation}
\frac12 \frac{d}{dt}\|(\omega e^{\omega_{\max} t}-\omega_{\min})_{-}\|_2^2\le  (\omega_{\max}-\omega,\omega e^{\omega_{\max}t}(\omega e^{\omega_{\max}, t}-\omega_{\min})_{-})\le 0, \label{minp22}
\end{equation}
where  we used \eqref{ole} and \eqref{oge0} to obtain the second inequality. Thus, by using the assumption \eqref{Asso}, we may conclude
\begin{equation}
\omega \ge \omega_{\min}e^{-\omega_{\max}t} \qquad \textrm{ a.e.  in } Q. \label{minpo}
\end{equation}
Therefore, \eqref{minmaxom} immediately follows.

\subsection{Proof of Lemma~\ref{n-ap}}\label{sub3.2}
In this subsection, we use both  $(\bc^m, \bd^m, \omega^m)$ and also $(\bv^m,b^m,\omega^{m})$ to denote a solution to $\mathcal{P}_{m,n,k}$ whose existence was established in  Lemma~\ref{m-ap}. Our goal is to let $m\to \infty$ to prove Lemma~\ref{n-ap}.
\subsubsection{Uniform $m$-independent estimates}
Repeating the same procedure as in Subsection~\ref{sub3.1}, we find that \begin{align}
\sup_{t\in(0,T)}\|\bv^m(t)\|_2^2 +\int_0^T
\int_{\Omega}T_k(\mu^{n,m})|\bD(\bv^m)|^2\; dx \; dt&\le
C, \label{apest1m}\\
\sup_{t\in (0,T)}|\bc^m(t)|&\le C,\label{apestsimm}\\
\sup_{t\in (0,T)}|\partial_{t} \bc^m(t)|&\le C(n,k).\label{apestsim2m}
\end{align}
Similarly, we obtain the identity
\begin{equation}
\begin{split}
(\partial_{t} b^m, b^m-b_{\Gamma}^n) - (b^m \bv^m, \nabla (b^m-b^n_{\Gamma}) + (T_n(\mu^{n,m})\nabla b^m, \nabla (b^m - b_{\Gamma}^n)\\
=-\left(b_+^{m} \omega^m + \frac{T_k(\mu^{n,m})|\bD(\bv^m)|^2}{1+n^{-1}|\bD(\bv^m)|^2}, b^m-b_{\Gamma}^n\right).
\end{split}\label{bm}
\end{equation}
By virtue of \eqref{apestsimm}, \eqref{minmaxom}, the smoothness of  $b_{\Gamma}^n$ and $b_0^{n,k}$, Young's inequality and Gronwall's lemma, we find that
\begin{equation}
\sup_{t\in (0,T)} \|b^m(t)\|_2^2 + \frac{1}{n} \int_0^T \|\nabla b^m\|_2^2 \; dt \le C(k,n). \label{aprbm}
\end{equation}
Then, it follows from \eqref{weak3-1l} that
\begin{equation}
\int_0^T \|\partial_{t} b^m\|_{W^{-1,2}_{\Gamma}(\Omega)}^2\; dt \le C(n,k).\label{timebm}
\end{equation}
Finally, setting $z:=\omega^m-\omega_{\Gamma}^m$ in \eqref{weak2-1l} we get the identity
\begin{equation}
\begin{split}
\frac12 &\frac{d}{dt} \| \omega^m-\omega_{\Gamma}^m\|_2^2  + \int_{\Omega} T_n(\mu^{n,m})|\nabla (\omega^m-\omega_{\Gamma}^m)|^2\; dx+((\omega^m)^2, \omega^m-\omega_{\Gamma}^m)\\
&= -(T_n(\mu^{n,m})\nabla \omega^m_{\Gamma}, \nabla (\omega^m - \omega_{\Gamma}^m)-(\partial_{t} \omega_{\Gamma}^m,\omega^m-\omega_{\Gamma}^m) -(\bv^m (\omega^m-\omega^m_{\Gamma}), \nabla  \omega^m_{\Gamma}). \label{ID1o2m}
\end{split}
\end{equation}
Then,  using \eqref{bd3}, \eqref{wboundar}, \eqref{minmaxom}, \eqref{dfmunm} and \eqref{apestsimm} we observe that
\begin{equation}
\frac{1}{n}\int_0^T \|\nabla \omega^m\|_2^2\; dt \le C(k,n), \label{nabom}
\end{equation}
and consequently we also have
\begin{equation}
\int_0^T \|\partial_{t} \omega^m\|_{W^{-1,2}_{\Gamma}(\Omega)}^2\; dt \le C(k,n).\label{timeom}
\end{equation}
\subsubsection{Taking the limit $m\to \infty$}
Using \eqref{wboundar}, \eqref{ole}, \eqref{minpo}, \eqref{apest1m}--\eqref{apestsim2m}, \eqref{aprbm}--\eqref{timebm} and \eqref{nabom}--\eqref{timeom} we can find a subsequence that we do not relabel such that
\begin{align}
\bc^m &\rightharpoonup^* \bc &&\textrm{weakly$^*$ in } W^{1,\infty}(0,T)^n,\label{cm1}\\
\bc^m &\to \bc &&\textrm{strongly in } \mathcal{C}(0,T)^n,\label{cm1s}\\
\bv^m&\to \bv &&\textrm{strongly in } \mathcal{C}(0,T; W^{1,2}_{\bn, \diver}), \label{cm3}\\
b^m-b_{\Gamma}^n&\rightharpoonup^* b-b_{\Gamma}^n &&\textrm{weakly$^*$ in } L^2(0,T; W^{1,2}_{\Gamma}(\Omega)) \cap L^{\infty}(0,T; L^2(\Omega)), \label{cm4}\\
\partial_{t} b^m &\rightharpoonup \partial_{t} b &&\textrm{weakly in } L^2(0,T; W^{-1,2}_{\Gamma}(\Omega)), \label{cm41}\\
b^m &\to b &&\textrm{strongly in } L^2(0,T; L^2(\Omega)), \label{cm42}\\
\omega^m-\omega_{\Gamma}^m &\rightharpoonup^* \omega-\omega_{\Gamma} &&\textrm{weakly$^*$ in } L^2(0,T; W^{1,2}_{\Gamma}(\Omega))\cap L^{\infty}(0,T;L^{\infty}(\Omega)),\label{cm5}\\
\partial_{t} \omega^m &\rightharpoonup \partial_{t} \omega &&\textrm{weakly in } L^2(0,T; W^{-1,2}_{\Gamma}(\Omega)), \label{cm5.1}\\
\omega^m &\to \omega &&\textrm{strongly in } L^2(0,T; L^2(\Omega)). \label{cm6}
\end{align}
Having these convergence results and the minimum principle \eqref{minpo} it is standard to let $m\to \infty$ in \eqref{weak1-1l}--\eqref{weak2-1l} to get \eqref{weak1-1m}--\eqref{dfmunn}, provided that we show the validity of \eqref{minbn}. Moreover, the attainment of the initial data \eqref{oinda2} can be deduced by standard tools.
\subsubsection{Minimum principle for $b$}
First, notice that letting $m\to \infty$ in \eqref{weak3-1l} without assuming the nonnegativity of $b$, we obtain
\begin{align}
&\left.
\begin{aligned}
&\langle \partial_{t} b, z\rangle -(b \bv, \nabla z) + (T_n(\tilde{\mu}^n) \nabla b, \nabla z) = (-b_+ \omega + T_k(\tilde{\mu}^n)|\bD(\bv)|^2, z)\\
&\qquad \textrm{ for all } z\in W^{1,2}_{\Gamma}(\Omega) \textrm{ and a.a. } t\in (0,T),
\end{aligned}\right. \label{weak3-1m*}
\end{align}
with $\tilde{\mu}^n$ is given as
\begin{equation}
\tilde{\mu}^n:=\frac{b_+}{\omega}+ \frac{1}{n}.\label{dfmunn*}
\end{equation}
Next, since $b_{\Gamma}^n\ge b_{\min}$ we can take $z:=b_{-}$ in \eqref{weak3-1m*}. First, the convective term vanishes and the third term on the left-hand side is nonnegative. Moreover, we see that the term on the right-hand side is non-positive and consequently by using \eqref{inn} we find that $b\ge 0$ almost everywhere in $Q$ and therefore we can replace $b_+$ by $b$. Next, setting $z:=e^{\omega_{\max}t} (b e^{\omega_{\max}t} - k^{-1})_{-}$ in \eqref{weak3-1m} (note that such a setting is admissible since $k\ge \frac{1}{b_{\min}}$) we derive an inequality (using the fact that $\omega \le \omega_{\max}$)
\begin{equation}
\frac{d}{dt}\|(b e^{t\omega_{\max}} - k^{-1})_{-}\|_2^2 \le 0,
\label{mpb}
\end{equation}
which implies after using \eqref{inn} the relation \eqref{minbn}.

\subsection{Proof of Lemma~\ref{k-lemma}}
In order to prove Lemma~\ref{k-lemma} we have to leave the standard $L^2$-theory. This is why we provide a more detailed proof in this part. Let $(\bv^n,b^n,\omega^n)$ be a solution to problem $\mathcal{P}_{n,k}$, whose existence is guaranteed by Lemma~\ref{n-ap}. Our goal is to let $n\to \infty$ in \eqref{weak1-1m}--\eqref{weak2-1m} to prove Lemma~\ref{k-lemma}.

\subsubsection{Uniform $n$-independent estimates}
In the same way as in the preceding subsections, we can obtain the energy identity
\begin{equation}
\|\bv^n(t)\|_2^2 + 2\int_0^t
\int_{\Omega}T_k(\mu^n)|\bD(\bv^n)|^2\; dx + (\bg^k(\cdot,b^n,\omega^n,\bv^n), \bv^n)_{\partial \Omega} \; d\tau=\|\bv^n_0\|_2^2. \label{Ee1}
\end{equation}
Hence, using \eqref{Assv},  \eqref{lowg} and \eqref{dfgk}, we
observe that
\begin{equation}
\begin{split}
&\sup_{t\in(0,T)}\|\bv^n(t)\|_2^2 + \int_Q T_k(\mu^n)|\bD(\bv^n)|^2\; dx \; dt\\
&\qquad +\int_0^T\int_{\partial \Omega} |\bg^k(\cdot,b^n,\omega^n,\bv^n) \cdot \bv^n|\; dS\; dt\le
C.
\end{split}\label{apest1n}
\end{equation}
Thus, using \eqref{minmaxon},  \eqref{minbn} and Korn's inequality \eqref{Korn}, we deduce that
\begin{equation}
\int_0^T \|\bv^n\|^2_{1,2} \; dt\le
C(k). \label{apest1n*}
\end{equation}
Next, having \eqref{apest1n} and \eqref{apest1n*} we deduce from
\eqref{weak1-1m} that (note that this estimate is valid because of the presence of the cut-off functions $G_k$ and $T_k$)
\begin{equation}
\int_0^T \|\partial_{t} \bv^n\|_{W^{-1,2}_{\bn,\diver}}^2 \; dt \le C(k).
\label{apestvt}
\end{equation}
Moreover, it follows from the standard interpolation inequality
\begin{equation}
\|u\|_{\frac{10}{3}}\le C\|u\|^{\frac{2}{5}}_2 \|u\|_{1,2}^{\frac{3}{5}}
\label{stintin}
\end{equation}
and from the estimates \eqref{apest1n} and \eqref{apest1n*} that
\begin{equation}
\int_0^T\|\bv^n\|_{\frac{10}{3}}^{\frac{10}{3}}\; dt \le C(k).
\label{apest1.1}
\end{equation}

Next, we focus on uniform estimates for $b^n$. First, for arbitrary $a>0$, we set $z:= T_a(b^n-b^n_{\Gamma})$ in \eqref{weak3-1m}. Using integration by parts, the fact that $\diver \bv^n =0$ and a simple algebraic manipulation, we find the identity
\begin{equation}
\begin{split}
&\left\langle \partial_{t} (b^n-b_{\Gamma}^n), T_a(b^n-b^n_{\Gamma}) \right \rangle + (b^n-b^n_{\Gamma},\omega^n T_a(b^n-b^n_{\Gamma}))\\
&\; + \int_{\Omega}T_n (\mu^n) |\nabla T_a(b^n-b^n_{\Gamma})|^2 \; dx =\left(-b^n_{\Gamma}\omega^n + \frac{T_k(\mu^n)|\bD(\bv^n)|^2}{1+n^{-1}|\bD(\bv^n)|^2}, T_a(b^n-b^n_{\Gamma})\right) \\
&\;- (T_n(\mu^n) \nabla b^n_{\Gamma}, \nabla T_a(b^n-b^n_{\Gamma})) - \left \langle \partial_{t} b^n_{\Gamma}, T_a(b^n-b^n_{\Gamma})\right \rangle - (\bv^n T_a(b^n-b^n_{\Gamma}),\nabla b^n_{\Gamma}).
\end{split}\label{ab1}
\end{equation}
First, we have that (recall the definition of $\Theta_a$ in \eqref{Theta})
$$
\left\langle \partial_{t} (b^n-b_{\Gamma}^n), T_a(b^n-b^n_{\Gamma}) \right \rangle= \frac{d}{dt} \|\Theta_a (b^n-b^n_{\Gamma})\|_1.
$$
Next, since $\omega^n \ge 0$ (see \eqref{minmaxon}) we observe that the second term in \eqref{ab1} is nonnegative. In addition, the integral with respect to time from the first term on the right-hand side of \eqref{ab1} can be estimated by using the assumption \eqref{bd4}$_2$,  the a~priori estimate \eqref{apest1n} and the maximum principle \eqref{minmaxon} as follows
$$
\int_0^t\left(-b^n_{\Gamma}\omega^n + \frac{T_k(\mu^n)|\bD(\bv^n)|^2}{1+n^{-1}|\bD(\bv^n)|^2}, T_a(b^n-b^n_{\Gamma})\right)\; d\tau \le a (\omega_{\max} b_{\max} |\Omega| T +C).
$$
For the last term on the right-hand side of \eqref{ab1}, we use \eqref{bboundar} and the uniform estimate \eqref{apest1n} to get (with the help of the H\"{o}lder inequality)
$$
-\int_0^t  (\bv^n T_a(b^n-b^n_{\Gamma}),\nabla b^n_{\Gamma})\; d\tau \le a\int_Q|\bv^n||\nabla b^n_{\Gamma}|\; dx\; dt \le Ca.
$$
Next, \eqref{bboundar} also implies that
$$
- \int_0^t \left \langle \partial_{\tau}b^n_{\Gamma}, T_a(b^n-b^n_{\Gamma})\right \rangle \; d\tau \le Ca.
$$
Finally, for the second term on the right-hand side of \eqref{ab1} we use Young's inequality and the definition of $T_a$ to conclude that
\begin{align*}
- (T_n(\mu^n) \nabla b^n_{\Gamma}, \nabla T_a(b^n-b^n_{\Gamma}))&\le \frac12 \int_{\Omega} T_n(\mu^n) |\nabla T_a(b^n-b^n_{\Gamma})|^2 \; dx \\
&\; + \frac12 \int_{|b^n-b^n_{\Gamma}|< a} T_n(\mu^n) |\nabla b^n_{\Gamma}|^2\; dx.
\end{align*}
Since
$$
\mu^n \le \frac{a+b_{\max}}{\omega_{\min}}+\frac{1}{n}
$$
on the set where $|b^n-b^n_{\Gamma}|\le a$, which follows from the definition of $\mu^n$ (see \eqref{dfmunn}), the minimum principle for $\omega$ (see \eqref{minmaxon}) and from the assumption \eqref{bd4}$_2$, we get  that
\begin{align*}
- (T_n(\mu^n) \nabla b^n_{\Gamma}, \nabla T_a(b^n-b^n_{\Gamma}))&\le \frac12 \int_{\Omega} T_n(\mu^n) |\nabla T_a(b^n-b^n_{\Gamma})|^2 \; dx \\
&\; + C(a,\omega_{\min},b_{\max}) \|\nabla b^n_{\Gamma}\|^2_2.
\end{align*}
Hence, inserting all of the above estimates into \eqref{ab1}, integrating with respect to time and using \eqref{bboundar}, we deduce that
\begin{equation}
\begin{split}
\sup_{t\in (0,T)} \|\Theta_a (b^n-b_{\Gamma}^n)\|_1 + \int_0^T \int_{\Omega} T_n(\mu^n) |\nabla T_a (b^n-b_{\Gamma}^n)|^2\; dx \; dt\\
\le C(a) + \|\Theta_a (b^{n,k}_0-b_{\Gamma}^n(0))\|_1. \label{startbn}
\end{split}
\end{equation}
Consequently, since \eqref{startbn} is valid for any $a\ge 0$, we find by using \eqref{bcon2}--\eqref{inn}, \eqref{bd4} and the properties of $T_a$ and $\Theta_a$ the following uniform estimate
\begin{equation}
\begin{split}
\sup_{t\in (0,T)} \|b^n\|_1 + \int_0^T \int_{\Omega} T_n(\mu^n) |\nabla T_a (b^n)|^2\; dx \; dt \le C(a) \qquad \textrm{ for all } a\ge 0.\label{startbn2}
\end{split}
\end{equation}
Thus, we control the gradient of $b^n$, uniformly with respect to $n$, on the set where $b^n$ is not large, say, for example on the set where $b^n\le b_{\max}$. To get also the control on the sets where $b^n\ge b_{\max}$, we set $z:=((b^n+1)^{-\lambda} - (b_{\max}+1)^{-\lambda})_{-}$ in \eqref{weak3-1m}, where $\lambda \in (0,1)$ is arbitrary. Such setting is possible because $|z|\le 1$ and also $z\equiv0$ if $b^n\le b_{\max}$. Therefore, denoting $\Omega_{-}:=\{x\in \Omega; z<0\}$ and using integration by parts to eliminate the convective term,  we obtain the following identity
\begin{equation}
\begin{split}
&\frac{d}{dt} \int_{\Omega_{-}} \frac{(b^n+1)^{1-\lambda}}{1-\lambda} - \frac{b^n+1}{(b_{\max}+1)^{\lambda}}- \frac{\lambda (b_{\max}+1)^{1-\lambda}}{1-\lambda}\; dx\\
 &\qquad - \lambda\int_{\Omega_{-}}T_n(\mu^n) \frac{|\nabla b^n|^2}{(b^n+1)^{1+\lambda}}\; dx \\
 &=\left (-b^n \omega^n + \frac{T_k(\mu^n)|\bD(\bv^n)|^2}{1+n^{-1}|\bD(\bv^n)|^2}, \left(\frac{1}{(b^n+1)^{\lambda}} - \frac{1}{(b_{\max}+1)^{\lambda}}\right)_{-}\right).
\end{split}
\label{nablab}
\end{equation}
Finally, we integrate the result with respect to time. First, the term on the right-hand side is uniformly bounded due to \eqref{startbn2}$_1$, \eqref{apest1n}, \eqref{minmaxon} and the fact that $|z|\le 1$. Also the first term on the left-hand side can be bounded after integration over time with help of \eqref{startbn2}$_1$. Consequently, we have a uniform ($n$ and $k$ independent) estimate
\begin{equation}
\begin{split}
\int_0^T \int_{\Omega_{-}}\frac{T_n(\mu^n)}{(1+b^n)^{1+\lambda}} |\nabla b^n|^2\; dx\; dt \le C(\lambda^{-1}) \textrm{ for all } \lambda \in (0,1).
\end{split}
\label{nablabfin*}
\end{equation}
Hence, setting $a:=b_{\max}$ in \eqref{startbn2} and adding the result to \eqref{nablabfin*} we find that
\begin{equation}
\begin{split}
\int_0^T \int_{\Omega}\frac{T_n(\mu^n)}{(1+b^n)^{1+\lambda}} |\nabla b^n|^2\; dx\; dt \le C(\lambda^{-1}) \textrm{ for all } \lambda \in (0,1).
\end{split}
\label{nablabfin}
\end{equation}
As we shall show later, this estimate is sufficient to take the limit $n\to \infty$ in \eqref{weak3-1m} because of the presence of $\frac{1}{k}$ in the estimate \eqref{minbn}.  However, it would not be sufficient to take the limit $k\to \infty$. Therefore, we improve the uniform estimate \eqref{nablabfin} so that it provides more information about the behavior for small $b^n$. Note that on the sets where $b^n\ge b_{\min}$, the uniform estimate of $\nabla b^n$ as well as $\bD(\bv^n)$ follows from \eqref{nablabfin} and \eqref{apest1n}, respectively. Note also that when deriving the estimate for $\nabla b^n$ on the sets where $b^n \le b_{\min}$, we obtain,  as a byproduct, the estimate for $\bD(\bv^n)$ that \emph{does not} depend on either $n$ or $k$. To do so, we set $z:=((b^n)^{-1} - (b_{\min})^{-1})_+$ in \eqref{weak3-1m}. Such a test function is well defined thanks to \eqref{minbn} and  $b^n=b_{\Gamma}^n\ge b_{\min}$ on $(0,T)\times \Gamma$ which follows from \eqref{bd4}. Thus, defining $\Omega_{\min}:=\{x\in \Omega; b^n \le b_{\min}\}$, we find (after using integration by parts to eliminate the convective term)  the following identity
\begin{equation}
\begin{split}
\frac{d}{dt} \int_{\Omega_{\min}}  \ln \left (\frac{b^n}{b_{\min}}\right) -\frac{b^n-b_{\min}}{b_{\min}}\; dx- \int_{\Omega_{\min}}\frac{T_n(\mu^n)}{(b^n)^2}|\nabla b^n|^2\; dx \\
 =\left (-b^n \omega^n + \frac{T_k(\mu^n)|\bD(\bv^n)|^2}{1+n^{-1}|\bD(\bv^n)|^2}, \left(\frac{1}{b^n} - \frac{1}{b_{\min}}\right)_{+}\right).
\end{split}
\label{nablabb}
\end{equation}
Hence, moving the terms with the corresponding signs to one side we get after integration with respect to time
\begin{equation}
\begin{split}
&\int_{\Omega} \left(\ln \left (\frac{b_{\min}}{b^n(t)}\right)\right)_{+}\; dx+ \int_0^t \int_{\Omega_{\min}}\frac{T_n(\mu^n)}{(b^n)^2}|\nabla b^n|^2 + \frac{T_k(\mu^n)|\bD(\bv^n)|^2}{b^n(1+n^{-1}|\bD(\bv^n)|^2)}\; dx\; d\tau  \\
 &=\int_{0}^t \int_{\Omega_{\min}} \omega^n \left(1-\frac{b^n}{b_{\min}}\right) + \frac{T_k(\mu^n)|\bD(\bv^n)|^2}{b_{\min}(1+n^{-1}|\bD(\bv^n)|^2)}\; dx \; d\tau\\
 &\qquad +\frac{1}{b_{\min}}\int_{\Omega}  (b_{\min}-b^n(t))_+ -(b_{\min}-b^{n,k}_0)_+\; +b_{\min} \left(\ln \left (\frac{b_{\min}}{b^{n,k}_0}\right)\right)_{+}\; dx.
\end{split}
\label{nablabb2}
\end{equation}
Consequently, using \eqref{apest1n},  \eqref{startbn2}$_1$, the maximum principle for $\omega^n$ \eqref{minmaxon} and
\eqref{nablabfin}, we deduce that, for all $\lambda \in (0,1]$,
\begin{equation}
\begin{split}
&\sup_{t\in (0,T)} \|\ln b^n(t)\,\|_1 + \int_0^T \int_{\Omega}\frac{T_n(\mu^n)}{(b^n)^{1+\lambda}}|\nabla b^n|^2 + \frac{T_k(\mu^n)|\bD(\bv^n)|^2}{b^n(1+n^{-1}|\bD(\bv^n)|^2)}\; dx\; d\tau  \\
 &\quad\le C(b_{\min},\omega_{\max}, \lambda^{-1}) +\|\ln b_0^{n,k}\,\|_1 \le C(\lambda^{-1}),
\end{split}
\label{nablafin2}
\end{equation}
where the last inequality follows from the assumption \eqref{bcon2}--\eqref{inn} and \eqref{Assk}.

It follows from \eqref{Ee1} (on the sets where $b^n\ge b_{\min}$) and the third term in \eqref{nablafin2} (on the sets where $b^n\le b_{\min}$) that
\begin{equation}\label{estpepa}
\int_0^T\it_{\Omega} \frac{|\bD(\bv^n)|^2}{1+n^{-1}|\bD(\bv^n)|^2 }\; dx\; d\tau \le C\,.
\end{equation}

Next, we focus on estimates for $T_n(\mu^n)$ that are uniform w.r.t. both $n$ and $k$. It directly follows from the definition \eqref{dfmunn} that
\begin{equation}
T_n\left(\frac{b^n}{\omega^n}\right)\le T_n(\mu^n)\le T_n\left(\frac{b^n}{\omega^n}\right)+ \frac{1}{n}.\label{tnh1}
\end{equation}
Since
\begin{equation}
\min \left\{1,\frac{1}{\omega^n}\right\} T_n (b^n) \le T_n\left(\frac{b^n}{\omega^n}\right)\le  \max \left(1,\frac{1}{\omega^n}\right) T_n(b^n),\label{tnh2}
\end{equation}
we observe, using \eqref{minmaxon}, \eqref{tnh1} and \eqref{tnh2}, that there exist positive constants  $C_1$ and $C_2$ independent of $k$ and $n$ such that
\begin{equation}
C_1 T_n (b^n) \le T_n\left(\mu^n \right)\le  C_2 T_n(b^n) + \frac{1}{n}.\label{tnh3}
\end{equation}
Thus, by virtue of \eqref{nablafin2}, we find that, for all $\lambda \in (0,1)$,
\begin{equation}
\begin{split}
\int_0^T \|\nabla (T_n(b_n))^{1-\frac{\lambda}{2}}\|_2^2\; dt &= C(\lambda)\int_0^T \int_{\Omega} \frac{|\nabla T_n(b^n)|^2}{(T_n(b^n))^{\lambda}}\;dx \; dt\\
&=C(\lambda)\int_0^T \int_{\Omega} \frac{T_n(b^n)}{(b^n)^{1+\lambda}}|\nabla T_n(b^n)|^2\;dx \; dt \le C(\lambda^{-1}). \label{estT}
\end{split}
\end{equation}
Combining this inequality with \eqref{startbn2}$_1$ (where we set $a=1$ for example), we obtain
\begin{equation}
\begin{split}
\sup_{t\in (0,T)} \|(T_n(b^n))^{1-\frac{\lambda}{2}}\|_1 + \int_0^T \|(T_n(b_n))^{1-\frac{\lambda}{2}}\|_{1,2}^2\; dt \le C(\lambda^{-1}). \label{estT2}
\end{split}
\end{equation}
Using the interpolation inequality
$$
 \|g\|_{\frac83}^{\frac83}\le C \|g\|_{1}^{\frac{2}{3}}  \|g\|_{1,2}^2,
$$
the equivalence \eqref{tnh3} and the fact that  \eqref{estT2} holds for all $\lambda\in (0,1)$, we finally conclude that
\begin{equation}
\int_0^T \|T_n(\mu^n)\|_{\frac{8-\lambda}{3}}^{\frac{8-\lambda}{3}} \; dt \le C(\lambda^{-1}) \qquad \textrm{ for all } \lambda\in (0,1).\label{apTn}
\end{equation}

We continue with  $k$-dependent estimates for $b^n$ that help us to establish the proper convergence results when $n\to \infty$. First, using the definition of $\mu^n$ (see \eqref{dfmunn}),  the minimum principle \eqref{minbn} and \eqref{nablafin2} we get
\begin{equation}
\int_0^T\int_{\Omega} \frac{|\nabla b^n|^2}{(b^n)^{1+\lambda}}\; dx \; dt \le C(\lambda^{-1},k), \label{bkn}
\end{equation}
which after using \eqref{startbn2}$_1$ leads to
\begin{equation}
\int_0^T \|(b^n)^{\frac{1-\lambda}{2}}\|_{1,2}^2  \; dt \le C(\lambda^{-1},k) \qquad \textrm{ for all } \lambda\in (0,1).\label{bkn2}
\end{equation}
Consequently, using  the embedding theorem we get
\begin{equation}
\int_0^T \|(b^n)^{1-\lambda}\|_{3}  \; dt \le C(\lambda^{-1},k) \qquad \textrm{ for all } \lambda\in (0,1).\label{bkn3}
\end{equation}
Finally, using the interpolation inequality
$$
\|(b^n)^{1-\lambda}\|^{\frac53}_{\frac53}\le \|(b^n)^{1-\lambda}\|_1^{\frac23} \|(b^n)^{1-\lambda}\|_3,
$$
the a~priori bound \eqref{startbn2}$_1$ and the estimate \eqref{bkn3}, we obtain
\begin{equation}
\int_0^T \|b^n\|_{\frac{5-\lambda}{3}}^{\frac{5-\lambda}{3}} \le C(\lambda^{-1},k) \qquad \textrm{ for all } \lambda\in (0,1).
\label{intb}
\end{equation}

Next, we focus on the estimates generated by the diffusion term in \eqref{weak3-1m}, see the second term in \eqref{nablafin2}. For any $q\in (1,2)$, we observe, by H\"{o}lder's inequality, that
\begin{equation}
\begin{split}
\label{eellm}
\int_Q |\sqrt{T_n(\mu^n)}&\nabla b^n|^q\; dx \; dt = \int_Q \left( \frac{T_n(\mu^n)|\nabla b^n|^2}{(b^n)^{1+\lambda}} \right)^{\frac{q}{2}} (b^n)^{\frac{(1+\lambda)q}{2}}\; dx \;dt \\
&\le \left( \int_Q \frac{T_n(\mu^n)|\nabla b^n|^2}{(b^n)^{1+\lambda}}\; dx \; dt \right)^{\frac{q}{2}} \left(\int_Q  |b^n|^{\frac{(1+\lambda)q}{2-q}}\; dx \;dt\right)^{\frac{2-q}{2}}\\
&\le C(\lambda^{-1},k),
\end{split}
\end{equation}
where the last inequality follows from \eqref{nablafin2} and \eqref{intb} provided that we are able to find $\lambda \in (0,1)$ such that
$$
\frac{(1+\lambda)q}{2-q} < \frac53.
$$
This is however possible whenever $q\in [1,5/4)$. Therefore, using this bound, we have
\begin{equation}
\begin{split}
\label{eellm2}
\int_Q |\sqrt{T_n(\mu^n)}\nabla b^n|^{\frac{5-\lambda}{4}}\; dx \; dt \le C(\lambda^{-1},k) \qquad \textrm{ for all } \lambda\in (0,1).
\end{split}
\end{equation}
Note that we conclude, as a direct consequence of \eqref{eellm2}, \eqref{minbn} and \eqref{intb}, that
\begin{equation}
\int_0^T \|b^n\|_{1,\frac{5-\lambda}{4}}^{\frac{5-\lambda}{4}} \; dt \le C(\lambda^{-1},k) \qquad \textrm{ for all } \lambda\in  (0,1).
\label{nablbnn}
\end{equation}
In addition,  using H\"{o}lder's inequality, \eqref{apTn} and \eqref{eellm2} we find that
\begin{equation}
\begin{split}
\label{eellm2*}
\int_Q |T_n(\mu^n)\nabla b^n|^{\frac{80-\lambda}{79}}\; dx \; dt \le C(\lambda^{-1},k) \qquad \textrm{ for all } \lambda\in (0,1).
\end{split}
\end{equation}
In order to obtain a uniform bound on $\partial_t b^n$, it remains to estimate the convective term $b^n\bv^n$. It however follows from \eqref{intb}
and \eqref{apest1.1}, by using H\"{o}lder's inequality, that
\begin{equation}
\begin{split}
\label{eellm2**}
\int_Q |b^n \bv^n|^{\frac{10-\lambda}{9}} \; dx \; dt \le C(\lambda^{-1},k)  \qquad \textrm{ for all } \lambda \in (0,1).
\end{split}
\end{equation}
Thus, using \eqref{weak3-1m} and the estimates \eqref{apest1n}, \eqref{eellm2*} and \eqref{eellm2**}, we conclude that
\begin{equation}
\begin{split}
\label{timebn}
\int_0^T  \|\partial_{t} b^n\|_{W_{\Gamma}^{-1,\frac{80-\lambda}{79}}}\; dt \le C(\lambda^{-1},k)  \qquad \textrm{ for all } \lambda \in (0,1).
\end{split}
\end{equation}

Finally, we derive uniform (yet $k$-dependent estimates) for $\omega^n$. Recall that the estimate \eqref{minmaxon} is uniform with respect to both $n$ and $k$. Hence, we focus on estimates for $\nabla \omega^n$. To do so, we set $u:=\omega^n - \omega_{\Gamma}$ in \eqref{weak2-1m}, and repeating the same procedure as in the proof of Lemma~\ref{m-ap} we get (see \eqref{ID1o2m})
\begin{equation}
\begin{split}
&\frac12 \frac{d}{dt} \| \omega^n-\omega_{\Gamma}\|_2^2  + \int_{\Omega} T_n(\mu^n)|\nabla (\omega^n-\omega_{\Gamma})|^2\; dx+((\omega^n)^2, \omega^n-\omega_{\Gamma})\\
&= -(T_n(\mu^n)\nabla \omega_{\Gamma}, \nabla (\omega^n - \omega_{\Gamma}))-(\partial_{t} \omega_{\Gamma},\omega^n-\omega_{\Gamma}) -(\bv^n (\omega^n-\omega_{\Gamma}), \nabla  \omega_{\Gamma}). \label{ID1o2n}
\end{split}
\end{equation}
Hence, using Young's inequality, \eqref{bd3} and \eqref{minmaxon} we get
\begin{equation}
\begin{split}
&\frac{d}{dt} \| \omega^n-\omega_{\Gamma}\|_2^2  + \int_{\Omega} T_n(\mu^n)|\nabla (\omega^n-\omega_{\Gamma})|^2\; dx\\
&\le C\left(1+\int_{\Omega} T_n(\mu^n)|\nabla \omega_{\Gamma}|^2 + |\partial_{t} \omega_{\Gamma}| + |\bv^n| |\nabla  \omega_{\Gamma}|\; dx\right)\\
&\le C(1+\|T_n(\mu^n)\|_{\frac{\beta_{\Gamma}}{\beta_{\Gamma}-2}}^{\frac{\beta_{\Gamma}}{\beta_{\Gamma}-2}}+ \|\omega_{\Gamma}\|_{1,\beta_{\Gamma}}^{\beta_{\Gamma}}+\|\partial_{t} \omega_{\Gamma}\|_1 +\|\bv^n\|_2^2 + \|\omega_{\Gamma}\|_{1,2}^2). \label{ID1o2n3}
\end{split}
\end{equation}
Since $\beta_{\Gamma}>\frac{16}{5}>2$ (see the assumption \eqref{bd1}), we have that
$$
\frac{\beta_{\Gamma}}{\beta_{\Gamma}-2}<\frac{8}{3}.
$$
Thus, integrating \eqref{ID1o2n3} with respect to time and using \eqref{Asso}, \eqref{bd1},  \eqref{bd3}, \eqref{apest1n} and \eqref{apTn} we deduce that
\begin{equation}
\begin{split}
 \int_Q T_n(\mu^n)|\nabla \omega^n|^2\; dx \; dt \le C, \label{apest2n}
\end{split}
\end{equation}
which,  after using \eqref{minbn},  implies
\begin{equation}
\begin{split}
\int_0^T \|\omega^n\|^2_{1,2}\; dt \le C(k). \label{apest2n*}
\end{split}
\end{equation}
In addition, we conclude from \eqref{apTn}, \eqref{apest2n} and H\"{o}lder's inequality that
\begin{equation}
\begin{split}
\int_0^T \|T_n(\mu^n)\nabla \omega^n\|_{\frac{16-\lambda}{11}}^{\frac{16-\lambda}{11}}\; dt \le C(\lambda^{-1}) \qquad \textrm{ for all } \lambda \in (0,1). \label{fluxon}
\end{split}
\end{equation}
Having the a~priori estimates \eqref{apest1n} and \eqref{fluxon} and the maximum principle \eqref{minmaxon}, we derive from \eqref{weak2-1m} the bound
\begin{equation}
\begin{split}
\int_0^T \|\partial_{t} (\omega^n)\|_{W^{-1,{\frac{16-\lambda}{11}}}_{\Gamma}}^{\frac{16-\lambda}{11}}\; dt \le C(\lambda^{-1}) \qquad \textrm{ for all } \lambda \in (0,1). \label{timeon}
\end{split}
\end{equation}

\subsubsection{On taking the limit $n\to \infty$}
Having \eqref{minmaxon}, \eqref{apest1n}--\eqref{apest1.1}, \eqref{nablbnn}, \eqref{timebn}, \eqref{apest2n}, \eqref{apest2n*} and \eqref{timeon}, we can let $n\to \infty$ and find a subsequence of $(\bv^n, b^n, \omega^n)$ that we do not relabel such that
\begin{align}
\bv^n &\rightharpoonup^* \bv &&\textrm{weakly$^*$ in } L^{\infty}(0,T;L^2_{\bn,\diver})\cap L^2(0,T;W^{1,2}_{\bn,\diver})
, \label{c1}\\
\partial_{t} \bv^n&\rightharpoonup \partial_{t} \bv &&\textrm{weakly in } L^2(0,T;W^{-1,2}_{\bn,\diver}), \label{c3}\\
b^n &\rightharpoonup b &&\textrm{weakly in } L^q(0,T;W^{1,q}(\Omega)) \textrm{ for all } q\in \left[1,\frac{5}{4}\right), \label{c4}\\
\partial_{t} b^n &\rightharpoonup \partial_{t} b &&\textrm{weakly in } \mathcal{M}(0,T;W^{-1,q}_{\Gamma}(\Omega)) \textrm{ for all } q\in \left[1, \frac{80}{79}\right), \label{c5}\\
\omega^n &\rightharpoonup^* \omega &&\textrm{weakly$^*$ in } L^2(0,T;W^{1,2}(\Omega))\cap L^{\infty}(0,T;L^{\infty}(\Omega)), \label{c2}\\
\partial_{t} \omega^n&\rightharpoonup \partial_{t} \omega &&\textrm{weakly in } L^q(0,T;W^{-1,q}_{\Gamma}(\Omega))\textrm{ for all } q\in \left[1, \frac{16}{11}\right). \label{c6}\\
\intertext{In addition, using the generalized version of Aubin-Lions lemma, we observe that there is  a subsequence (that is again not relabelled) such that, for any $\alpha \in (0,1)$,}
b^n &\to b &&\textrm{strongly in } L^{\frac{5}{4}}(0,T;W^{\alpha,\frac{5}{4}}(\Omega)), \label{sc1}\\
\omega^n &\to \omega &&\textrm{strongly  in } L^2(0,T;W^{\alpha,2}(\Omega)), \label{sc2}\\
\bv^n &\to \bv &&\textrm{strongly in  } L^2(0,T;W^{\alpha,2}(\Omega)^3 \cap W^{1,2}_{\bn,\diver}).\label{sc3}\\
\intertext{Thus, using the trace theorem \eqref{trace}, we find that}
\bv^n &\to \bv &&\textrm{strongly in } L^{2}(0,T;L^{2}(\partial \Omega)^3). \label{bc1*}\\
b^n &\to b &&\textrm{strongly in } L^{\frac{5}{4}}(0,T;L^{\frac{5}{4}}(\partial \Omega)). \label{bc2*}\\
\omega^n &\to \omega &&\textrm{strongly in } L^{2}(0,T;L^{2}(\partial \Omega)). \label{bc3*}\\
\intertext{Moreover, using \eqref{sc1}--\eqref{sc3}, there is a subsequence of $(\bv^n, b^n, \omega^n)$ (again not relabelled) such that}
\bv^n &\rightarrow \bv &&\textrm{a.e. in  }Q \quad \textrm{and a.e. in  }(0,T)\times \partial \Omega, \label{pc1}\\
\omega^n &\rightarrow \omega &&\textrm{a.e. in  }Q\quad \textrm{and a.e. in  }(0,T)\times \partial \Omega, \label{pc2}\\
b^n &\rightarrow b &&\textrm{a.e. in  }Q\quad \textrm{and a.e. in  }(0,T)\times \partial \Omega. \label{pc3}
\end{align}
Thus, having \eqref{c1}, \eqref{c3}, \eqref{sc3}, \eqref{bc1*} and \eqref{pc1}--\eqref{pc3}, and using the continuity of the cut-off functions $T_k, G_k$ and also the fact that $\bg^k$ is bounded and continuous with respect to $(\bv,b,\omega)$, it is easy to let $n\to \infty$ in \eqref{weak1-1m} and obtain \eqref{weak1-1k}. Moreover, at this level of approximation, it  is standard to show the attainment of $\bv_0$, i.e., to prove \eqref{oinda333}$_1$.

Next, in order to identify the limit of \eqref{weak2-1m} as $n\to \infty$, we notice that \eqref{fluxon} implies that (for a subsequence)
\begin{align}
T_n(\mu^n)\nabla \omega^n &\rightharpoonup \overline{\mu \nabla \omega} &&\textrm{weakly in } L^q(0,T; W^{1,q}(\Omega)) \quad \textrm{ for all } q\in \left[1,\frac{16}{11}\right). \label{co1}
\end{align}
Thus, using in addition \eqref{c6}, \eqref{sc2}, \eqref{sc3} and \eqref{minmaxon}, and letting $n\to \infty$ in \eqref{weak2-1m} we obtain
\begin{equation}\label{weak2-1m*}
\begin{split}
&\langle \partial_{t} \omega, z\rangle -(\bv\omega, \nabla z) + (\overline{\mu \nabla \omega}, \nabla z) = -(\omega^2,z)\\
&\qquad \textrm{ for all } z\in W^{1,\infty}_{\Gamma}(\Omega) \textrm{ and a.a. } t\in (0,T).
\end{split}
\end{equation}
Thus, \eqref{weak2-1m*} leads to \eqref{weak2-1k} once we show that
\begin{equation}\label{dream}
\overline{\mu \nabla \omega}= \mu \nabla \omega \qquad \textrm{ a.e. in } Q.
\end{equation}
To do so, we use \eqref{apTn}, \eqref{minmaxon} and \eqref{pc2}--\eqref{pc3} to conclude that
\begin{align}
T_n(\mu_n) &\to \mu &&\textrm{strongly in } L^q(0,T; L^q(\Omega)) \quad \textrm{ for all } q\in \left[1,\frac{8}{3}\right).\label{co2}
\end{align}
Then, \eqref{dream} is a direct consequence of \eqref{c2} and \eqref{co2}. Hence \eqref{weak2-1m*} is nothing else than \eqref{weak2-1k}. Finally, the attainment of the initial condition for $\omega$, see \eqref{oinda333}$_2$, can be proven in a standard way.

Finally, we focus on obtaining the limit as $n\to \infty$ in \eqref{weak3-1m}. First, we identify the weak limit of the diffusion term. It follows from \eqref{pc1}, \eqref{co2} and Vitali's lemma that, for all $\alpha\in (0,1)$,
\begin{align}
(T_n(\mu_n))^{\alpha} &\to \mu^{\alpha} &&\textrm{strongly in } L^{\frac{8}{3}}(0,T; L^{\frac{8}{3}}(\Omega)).\label{co2*}
\end{align}
In addition, combining \eqref{nablbnn} and \eqref{eellm2*} we obtain for all $\alpha\in [0,1]$ that\footnote{We use the fact that $\frac{81}{80} < \frac{80}{79}$.}
\begin{align}
(T_n(\mu^n))^{\alpha} \nabla b^n &\rightharpoonup \overline{\mu^{\alpha} \nabla b} &&\textrm{weakly in  } L^{\frac{81}{80}}(0,T;L^{\frac{81}{80}}(\Omega)). \label{bdis}
\end{align}
Our goal is  identify for all $\alpha \in [0,1]$
\begin{equation}
\overline{\mu^{\alpha} \nabla b}=\mu^{\alpha} \nabla b \qquad \textrm{ a.e. in } Q. \label{dream2}
\end{equation}
Note that in order to take the limit in \eqref{weak3-1m} it is enough to show \eqref{dream2} only for $\alpha=1$. To prove it, we proceed inductively. We define $h:=\frac{1}{81}$, $\alpha_0=0$ and $\alpha_{i+1}=\alpha_i + h$. Note that \eqref{dream2} holds for $\alpha_0$ and we want to show that if it holds for $\alpha_i$ it also holds for $\alpha_{i+1}$ for all $i=1,\ldots, 81$. Thus, assume that \eqref{dream2} holds for $\alpha_i$. Then
$$
\overline{\mu^{\alpha_{i+1}} \nabla b}=\overline{\mu^h \mu^{\alpha_i} \nabla b}\overset{\eqref{co2*},\eqref{bdis}}=\mu^h \overline{\mu^{\alpha_i} \nabla b}\overset{\eqref{dream2}}=\mu^h \mu^{\alpha_i} \nabla b=\mu^{\alpha_{i+1}}\nabla b.
$$
Hence, setting $i=81$, we get \eqref{dream2} with $\alpha=1$. Thus, using \eqref{eellm2*} and \eqref{bdis} we finally, observe that
\begin{align}
T_n(\mu_n)\nabla b &\rightharpoonup \mu \nabla b &&\textrm{weakly in  } L^q(0,T;L^q(\Omega)^3) \textrm{ for all } q\in \left[1, \frac{80}{79}\right).\label{c9}
\end{align}
Next, we focus on the convergence properties of the second term on the right-hand side of \eqref{weak3-1m}. First, setting $\bw:=\bv$ in \eqref{weak1-1k} and integrating the result with respect to time over $(0,T)$ we get the energy identity
\begin{equation}
\|\bv(T)\|_2^2 + 2\int_0^T
\int_{\Omega}T_k(\mu)|\bD(\bv)|^2\; dx + (\bg^k(\cdot, b,\omega, \bv), \bv)_{\partial \Omega} \; dt=\|\bv_0\|_2^2. \label{Ee1k}
\end{equation}
Hence, setting $t:=T$ in \eqref{Ee1} and letting $n\to \infty$ we get by using \eqref{pc1}--\eqref{pc3}, the boundedness of $\bg^k$, the weak-lower semicontinuity of norms and \eqref{Ee1k} that
\begin{equation}
\limsup_{n\to \infty}\int_0^T
\int_{\Omega}T_k(\mu^n)|\bD(\bv^n)|^2\; dx\; dt\le\int_0^T
\int_{\Omega}T_k(\mu)|\bD(\bv)|^2\; dx\; dt. \label{limsup}
\end{equation}
Since it directly follows from \eqref{c1} and \eqref{pc2}--\eqref{pc3} that
\begin{align}
\sqrt{T_k(\mu^n)}\bD(\bv^n)&\rightharpoonup \sqrt{T_k(\mu)}\bD(\bv) &&\textrm{weakly in } L^2(0,T;L^{2}(\Omega)^{3\times 3}), \label{c11}\\
\intertext{we immediately observe from \eqref{limsup} and \eqref{c11}, by referring to lower-semincontinuity of the $L^2$-norm, that}
T_k(\mu^n)|\bD(\bv^n)|^2&\to T_k(\mu)|\bD(\bv)|^2 &&\textrm{strongly in } L^1(0,T;L^{1}(\Omega)).\label{strongD}
\end{align}
Moreover, having \eqref{strongD} we can go back to \eqref{weak3-1m} and strengthen the convergence result \eqref{c5} to the following one
\begin{align}
\partial_{t} b^n &\rightharpoonup \partial_{t} b &&\textrm{weakly in } L^1(0,T;W^{-1,q}_{\Gamma}(\Omega)) \textrm{ for all } q\in \left[1, \frac{80}{79}\right). \label{c5bet}
\end{align}
Thus, having \eqref{c9}, \eqref{strongD} and \eqref{c5bet} it is easy to let $n\to \infty$ in \eqref{weak2-1m} to get \eqref{weak2-1k}. Moreover, one can deduce \eqref{oinda333}$_3$, but we postpone the proof of this to the proof of the main theorem where an even more difficult case is treated. In addition, one can use weak lower semicontinuity and Fatou's lemma and let $n\to \infty$ in \eqref{apest1n}, \eqref{startbn2}, \eqref{nablabfin}, \eqref{nablafin2}, \eqref{estpepa}, \eqref{apTn} and  \eqref{apest2n}  to get \eqref{finalap}.

\section{Proof of the main theorem}
\label{SMT}
For arbitrary $k\ge \frac{1}{b_{\min}}$, let us denote by $(\bv^k,b^k,\omega^k)$  the solution to the problem $\mathcal{P}_k$,   whose existence and uniform estimates are established in Lemma~\ref{k-lemma}. Next, we investigate the behavior of such solutions when $k \to \infty$ in order to prove the main result of the paper, i.e.,  Theorem~\ref{TH1}.

\subsection{Reconstruction of the pressure}
We start with the reconstruction of the pressure $p^k$. Note, that for the proof of Lemma~\ref{k-lemma} we required that $\Omega$ is only a Lipschitz domain. However, once we need to have control of the pressure, we use $W^{2,p}$ regularity for the Laplace operation (see \eqref{nabla2}) for which $\Omega\in \mathcal{C}^{1,1}$ is required (see also \cite{BuMaRa09,BuGwMaSw12} for further discussion how the regularity of $\Omega$ can be in certain cases weakened). Thus, following \cite{BuGwMaSw12}, we can show the existence of the pressure $p^k\in L^{2}(0,T;L_0^2(\Omega))$  such that, for almost all $t\in (0,T)$,
\begin{equation}
\begin{aligned}
&\langle \partial_{t} \bv^k, \bw \rangle
-\left(G_k(|\bv^k|^2)\bv^k\otimes \bv^k,
\nabla \bw\right) + (\bg^k(\cdot,b^k, \omega^k,\bv^k),
\bw)_{\partial \Omega} \\
&\qquad +\left(T_k(\mu^k)\bD(\bv^k),\bD(\bw)\right) =(p^k,\diver \bw)
\qquad \textrm{ for all }
\bw \in W^{1,2}_{\bn}.
\end{aligned}\label{prk}
\end{equation}
Moreover, the pressure can be decomposed as
\begin{equation}
p^k=p_1^k+p_2^k+p_3^k \qquad \textrm{such that }
\int_{\Omega} p_1^k \; dx=\int_{\Omega} p_2^k \; dx=\int_{\Omega} p_3^k \; dx=0 \label{decomp}
\end{equation}
and, for almost all time $t\in (0,T)$, $p_1^k$, $p_2^k$ and $p_3^k$ satisfy
\begin{equation}
\begin{aligned}
(p_1^k,\triangle \varphi)&=\left(T_k(\mu^k)\bD(\bv^k),\bD(\nabla \varphi)\right) &&\textrm{for all } \varphi\in W^{2,2}(\Omega); \;\nabla \varphi \in W^{1,2}_{\bn},\\
(p_2^k,\triangle \varphi)&=-\left(G_k(|\bv^k|^2)\bv^k\otimes \bv^k,\nabla^2 \varphi \right) &&\textrm{for all } \varphi\in W^{2,2}(\Omega);\;\nabla \varphi \in W^{1,2}_{\bn},\\
(p_3^k,\triangle \varphi)&=(\bg^k(\cdot,b^k,\omega^k,\bv^k),
\nabla \varphi)_{\partial \Omega} &&\textrm{for all } \varphi\in W^{2,2}(\Omega);\;\nabla \varphi \in W^{1,2}_{\bn}.
\end{aligned}\label{prk2}
\end{equation}

\subsection{Uniform estimates}
First, we recall \eqref{finalap} and derive further $k$-independent a~priori estimates. Note first that it follows from \eqref{finalap}, \eqref{estpepa} and \eqref{minmaxon} that
\begin{equation}
\begin{split}
\sup_{t\in (0,T)} \|\bv^k(t)\|_2^2 + \int_0^T \int_{\Omega}(1+T_k(\mu^k))|\bD(\bv^k)|^2\; dx \; dt\\
+\int_0^T\int_{\partial \Omega} |\bg^k(\cdot,\bv^k)\cdot \bv^k|\; dS \; dt \le C. \label{ae1}
\end{split}
\end{equation}
Using also the trace theorem (see \eqref{trace}) and the embedding theorem, we have
\begin{equation}\label{traceint}
\|\bv^k\|_{L^{\frac{8}{3}}(\partial \Omega)^3}^{\frac{8}{3}}\le C \|\bv^k\|_{W^{\frac{1}{4},2}(\partial \Omega)^3}^{\frac{8}{3}}\le C \|\bv^k\|_{W^{\frac{3}{4},2}(\Omega)^3}^{\frac{8}{3}}\le C \|\bv^k\|_{W^{1,2}(\Omega)^3}^{2} \|\bv^k\|_{L^{2}(\Omega)^3}^{\frac{2}{3}}.
\end{equation}
Hence, it follows from \eqref{ae1}, Korn's inequality \eqref{Korn} and the interpolation inequality \eqref{stintin} that
\begin{equation}
\int_0^T \|\bv^k\|_{1,2}^2 + \|\bv^k\|_{\frac{10}{3}}^{\frac{10}{3}} + \|\bv^k\|_{L^{\frac{8}{3}}(\partial \Omega)^3}^{\frac{8}{3}}\; dt \le C. \label{ae2}
\end{equation}
Consequently, by virtue of \eqref{upg} and \eqref{dfgk}, we conclude from \eqref{ae1} and \eqref{ae2} that
\begin{equation}
\int_0^T \|\bg^k(\cdot, b^k,\omega^k,\bv^k)\|^{\beta_{g}}_{L^{\beta_{g}}(\partial \Omega)^3} \; dt \le C. \label{ae3}
\end{equation}
Further, it directly follows from \eqref{finalap} and \eqref{minmaxok} that
\begin{equation}
\int_0^T \|T_k(\mu^k)\|_{\frac{8-\lambda}{3}}^{\frac{8-\lambda}{3}} \; dt \le C(\lambda^{-1}) \qquad \textrm{ for all } \lambda \in (0,1)\label{ae4}
\end{equation}
and then by using H\"{o}lder's inequality and \eqref{ae1} we find that
\begin{equation}
\int_0^T \|T_k(\mu^k)\bD(\bv^k)\|_{\frac{16-\lambda}{11}}^{\frac{16-\lambda}{11}} \; dt \le C(\lambda^{-1}) \qquad \textrm{ for all } \lambda \in (0,1). \label{ae5}
\end{equation}

Next, we focus on uniform estimates for the pressure. We derive only  the estimates for $p_3^k$ since the procedure developed in \cite{BuMaRa09} can be taken step by step to obtain the estimates for $p_k^1$ and $p_k^2$, which we provide here without proofs. Recalling that $\Omega \in \mathcal{C}^{1,1}$ we apply the $L^q$-theory for the Poisson equation (see \eqref{laplace}--\eqref{nabla2}). Consequently, we find a function $\varphi$ solving, for almost all $t\in (0,T)$,
\begin{align*}
\triangle \varphi &= |p_3^k|^{\frac{3\beta_{g}}{2}-2}p_3^k-\frac{1}{|\Omega|}\int_{\Omega}|p_3^k|^{\frac{3\beta_{g}}{2}-2}p_3^k\; dx &&\textrm{ in }\Omega,\\
\nabla \varphi \cdot \bn &= 0 &&\textrm{ on }\partial \Omega,\\
\int_{\Omega} \varphi \; dx &=0,
\end{align*}
and satisfying the following estimate
\begin{equation}
\|\varphi\|_{2,\frac{3\beta_{g}}{3\beta_{g}-2}}\le C \|p_3^k\|_{\frac{3\beta_{g}}{2}}^{\frac{3\beta_{g}}{2}-1}.\label{varp}
\end{equation}
Thus, using such a $\varphi$ in \eqref{prk2}$_3$ we get the identity (we use the fact that $\int_{\Omega} p_3^k = 0$)
\begin{equation*}
\|p_3^k\|_{\frac{3\beta_{g}}{2}}^{\frac{3\beta_{g}}{2}}=(\bg^k(\cdot,b^k,\omega^k, \bv^k),\nabla \varphi)_{\partial \Omega}.
\end{equation*}
Hence, applying H\"{o}lder's inequality, \eqref{varp}, the trace theorem \eqref{trace} and the embedding theorem, we obtain
\begin{equation}
\begin{split}
\|p_3^k\|_{\frac{3\beta_{g}}{2}}^{\frac{3\beta_{g}}{2}}&\le \|\bg^k\|_{L^{\beta_{g}}(\partial \Omega)^3} \|\nabla \varphi\|_{L^{\frac{\beta_{g}}{\beta_{g}-1}}(\partial \Omega)^3}\\
&\le C\|\bg^k\|_{L^{\beta_{g}}(\partial \Omega)^3} \|\varphi \|_{2,\frac{3\beta_{g}}{3\beta_{g}-2}}\le C\|\bg^k\|_{L^{\beta_{g}}(\partial \Omega)^3}\|p_3^k\|_{\frac{3\beta_{g}}{2}}^{\frac{3\beta_{g}}{2}-1}.
\end{split}\label{ae9*}
\end{equation}
By a simple algebraic manipulation, it follows from \eqref{ae9*} and \eqref{ae3} that
\begin{equation}
\begin{split}
\int_0^T\|p_3^k\|_{\frac{3\beta_{g}}{2}}^{\beta_{g}}\; dt &\le C.
\end{split}\label{ae9}
\end{equation}
In addition, going back to \eqref{prk2}$_3$, we see that for almost all times $t\in (0,T)$, the pressure $p_3^k$ is a harmonic function in $\Omega$. Thus, it directly follows from \eqref{ae9} that
\begin{equation}
\begin{split}
\int_0^T\|p_3^k\|_{L^{\infty}(\Omega')}^{\beta_{g}}\; dt &\le C(\Omega') \quad \textrm{ for all } \Omega' \subset \overline{\Omega'} \subset \Omega.
\end{split}\label{ae9**}
\end{equation}
Similarly,  \eqref{prk2}, \eqref{ae2}, \eqref{ae5} and the fact that $|G_k|\le 1$ for all $k\in \mathbb{N}$ imply
\begin{equation}
\int_0^T \|p_1^k\|_{\frac{16-\lambda}{11}}^{\frac{16-\lambda}{11}} +\|p_2^k\|_{\frac{5}{3}}^{\frac{5}{3}}\; dt \le C(\lambda^{-1}) \qquad \textrm{ for all } \lambda \in (0,1). \label{ae6}
\end{equation}
Finally, recalling the definition of $q_{\min}:=\min\{\beta_{\bg}, 16/11\}$ (see the formulation of Theorem \ref{TH1}), using \eqref{ae2}--\eqref{ae3}, \eqref{ae5}, \eqref{ae9} and \eqref{ae6}, we can derive from the identity  \eqref{prk}, \eqref{decomp} that for all $q\in [1,q_{\min})$ there holds
\begin{equation}
\int_0^T \|\partial_{t} \bv^k\|_{W^{-1,q}_{\bn}}^{q}\; dt \le C(q).\label{ae10}
\end{equation}

We continue with establishing the estimates for $b^k$. Having \eqref{finalap} and \eqref{minmaxok} it is easy to deduce with the help of H\"{o}lder's inequality that
\begin{equation}
\begin{split}
\sup_{t\in (0,T)} \left( \|b^k(t)\|_1 + \|\ln b^k(t)\|_1\right)+ \int_0^T \|\nabla b^k\|_{2-\lambda}^{2-\lambda} + \|b^k\|_{\frac{8-\lambda}{3}}^{\frac{8-\lambda}{3}} \; dt\\
+ \int_0^T \|\mu^k \nabla b^k\|_{\frac{8-\lambda}{7}}^{\frac{8-\lambda}{7}} \; dt\le C(\lambda^{-1})\qquad \textrm{ for all } \lambda \in (0,1).\label{ae11}
\end{split}
\end{equation}
Finally, we deduce from \eqref{weak3-1k}, using the estimates \eqref{ae1}, \eqref{ae2}, \eqref{ae11} and the maximum principle \eqref{minmaxok}, that
\begin{equation}
\begin{split}
\int_0^T \|\partial_{t} b^k\|_{W^{-1,\frac{8-\lambda}{7}}_{\Gamma}}\; dt \le C(\lambda^{-1})\qquad \textrm{ for all } \lambda \in (0,1).\label{ae12}
\end{split}
\end{equation}

We end this subsection by proving the uniform estimates for $\omega^k$. Since $b^k$ can vanish on a set of zero measure, we are not be able to bound $\omega^k$ in a Sobolev space and consequently to identify the trace of $\omega^k$ when $k\to \infty$. On the other hand, we show that $b^k\omega^k$ has such desired properties, which will be sufficient for identifying the limit as $k\to \infty$. Note first that \eqref{finalap} and \eqref{minmaxok} imply that
\begin{equation}
\sup_{t\in (0,T)}\|\omega^k(t)\|_{\infty} + \int_0^T \int_{\Omega} b^k|\nabla \omega^k|^2\; dx \; dt \le C. \label{ae13}
\end{equation}
Hence, using \eqref{ae11}, \eqref{ae13} and H\"{o}lder's inequality we get
\begin{equation}
\int_0^T \|b^k\nabla \omega^k\|^{\frac{16-\lambda}{11}}_{\frac{16-\lambda}{11}} \; dt \le C(\lambda^{-1}) \qquad \textrm{ for all } \lambda \in (0,1). \label{ae14}
\end{equation}
Next, using \eqref{ae11}, \eqref{ae14} and \eqref{minmaxok} we deduce that
\begin{equation}
\int_0^T \|\nabla (b^k \omega^k)\|^{\frac{16-\lambda}{11}}_{\frac{16-\lambda}{11}} \; dt \le C(\lambda^{-1}) \qquad \textrm{ for all } \lambda \in (0,1) \label{ae15}
\end{equation}
and also
\begin{equation}
\int_0^T \left\|\nabla \left(\frac{b^k \omega^k}{b^k+1}\right)\right\|^{\frac{16-\lambda}{11}}_{\frac{16-\lambda}{11}} \; dt \le C(\lambda^{-1}) \qquad \textrm{ for all } \lambda \in (0,1). \label{ae16}
\end{equation}
Moreover, taking \eqref{ae2}, \eqref{minmaxok} and \eqref{ae14} into account we conclude from \eqref{weak2-1k} that
\begin{equation}
\int_0^T \|\partial_{t} \omega^k\|_{W^{-1,\frac{16-\lambda}{11}}_{\Gamma}}^{\frac{16-\lambda}{11}}\; dt \le C(\lambda^{-1}) \qquad \textrm{ for all } \lambda\in (0,1).\label{ae17}
\end{equation}
In addition, since we shall need to identify the trace of $\omega$ on $\Gamma$ with $\omega_{\Gamma}$, we can in fact extend the estimate \eqref{ae15} to
\begin{equation}
\int_0^T \|b^k (\omega^k-\omega_{\Gamma})\|^{\frac{16-\lambda}{11}}_{W^{1,\frac{16-\lambda}{11}}_{\Gamma}} \; dt \le C(\lambda^{-1}) \qquad \textrm{ for all } \lambda \in (0,1), \label{ae15**}
\end{equation}
where we have used the assumption $\beta_{\Gamma}>\frac{16}{5}$ on $\omega_{\Gamma}$ (see \eqref{bd1}) and the uniform estimate \eqref{ae11} for $b^k$.

\subsection{On taking the limit $k\to \infty$}
In this final subsection we let $k\to \infty$ and show that the limit object solves the original problem. Using the uniform estimates \eqref{ae1}, \eqref{ae3}, \eqref{ae5}--\eqref{ae15**} and \eqref{minmaxok}, we can extract a subsequence that we do not relabel and  find a quintuple $(\bv,b,\omega,p,\bs$) with $p=p_1+p_2+p_3$ such that
\begin{align}
\bv^k &\rightharpoonup^* \bv &&\textrm{weakly$^*$ in } L^{\infty}(0,T;L^2_{\bn,\diver})\cap L^2(0,T;W^{1,2}_{\bn,\diver})
, \label{c1k}\\
\partial_{t} \bv^k&\rightharpoonup \partial_{t} \bv &&\textrm{weakly in } L^{q}(0,T;W^{-1,q}_{\bn})\textrm{ for all } q\in [1,q_{\min}), \label{c3k}\\
b^k-b_{\Gamma} &\rightharpoonup b-b_{\Gamma} &&\textrm{weakly in } L^q(0,T;W^{1,q}_{\Gamma}(\Omega)) \textrm{ for all } q\in \left[1,2\right), \label{c4k}\\
\partial_{t} b^k &\rightharpoonup \partial_{t} b &&\textrm{weakly in } \mathcal{M}(0,T;W^{-1,q}_{\Gamma}(\Omega)) \textrm{ for all } q\in \left[1, \frac{8}{7}\right), \label{c5k}\\
\omega^k &\rightharpoonup^* \omega &&\textrm{weakly$^*$ in } L^{\infty}(0,T;L^{\infty}(\Omega)), \label{c2k}\\
\partial_{t} \omega^k&\rightharpoonup \partial_{t} \omega &&\textrm{weakly in } L^q(0,T;W^{-1,q}_{\Gamma}(\Omega))\textrm{ for all } q\in \left[1, \frac{16}{11}\right), \label{c6k}\\
\bv^k &\rightharpoonup \bv &&\textrm{weakly in } L^{\frac{10}{3}}(0,T;L^{\frac{10}{3}}(\Omega)^3), \label{c3k*}\\
b^k &\rightharpoonup b &&\textrm{weakly in } L^q(0,T;L^{q}(\Omega)) \textrm{ for all } q\in \left[1,\frac{8}{3}\right). \label{c4k*}\\
p^k_{1}&\rightharpoonup p_{1} &&\textrm{weakly in } L^q(0,T;L^q(\Omega))\textrm{ for all } q\in \left[1, \frac{16}{11}\right). \label{c6kp}\\
p^k_2 &\rightharpoonup p_2 &&\textrm{weakly in } L^{\frac{5}{3}}(0,T;L^{\frac{5}{3}}(\Omega)), \label{c3k*p}\\
p^k_3 &\rightharpoonup p_3 &&\textrm{weakly in } L^{\beta_{g}}(0,T;L^{\frac{3\beta_{g}}{2}}(\Omega)), \label{c4k*p}\\
p^k_3 &\rightharpoonup p_3 &&\textrm{weakly in } L^{\beta_{g}}(0,T;L^{\infty}(\Omega'))\textrm{ for all } \overline{\Omega'}\subset \Omega, \label{c4k*p*}\\
\bg^k &\rightharpoonup \bs &&\textrm{weakly in } L^{\beta_{\bg}}(0,T;L^{\beta_{\bg}}(\partial\Omega)^{3}). \label{gconc}
\end{align}
In addition, using the generalized version of the Aubin-Lions lemma, we observe that there is  a subsequence (that is again not relabelled) such that, for all $\alpha\in (0,1)$ and all $q\in [1,2)$,
\begin{align}
b^k &\to b &&\textrm{strongly in } L^{q}(0,T;W^{\alpha,q}(\Omega)), \label{sc1k***}\\
\bv^k &\to \bv &&\textrm{strongly in  } L^{2}(0,T;W^{\alpha,2}(\Omega)^3 \cap W^{1,2}_{\bn,\diver}).\label{sc3k***}
\end{align}
Thus, combining these convergence results with the a~priori estimates \eqref{ae1} and \eqref{ae11} and using the trace theorem (see \eqref{trace}) we obtain (again for a subsequence)
\begin{align}
b^k &\to b \quad \, \textrm{ a.e. in } Q \textrm{ and strongly in } L^{q}(0,T;L^{q}(\Omega)) \textrm{ for all } q\in \left[1,\frac{8}{3}\right), \label{sc1k}\\
\bv^k &\to \bv \quad \textrm{ a.e. in } Q \textrm{ and strongly in } L^{q}(0,T;L^{q}(\Omega)^3) \textrm{ for all } q\in \left[1,\frac{10}{3}\right).\label{sc3k}\\
b^k &\to b \quad \, \textrm{ a.e. in } (0,T)\times \partial \Omega \textrm{ and strongly in } L^{1}(0,T;L^{1}(\partial\Omega)), \label{sc1ktr}\\
\bv^k &\to \bv \quad \textrm{ a.e. in } (0,T)\times \partial \Omega \textrm{ and strongly in } L^{1}(0,T;L^{1}(\partial\Omega)^3).\label{sc3ktr}
\end{align}
Consequently, it follows from \eqref{sc1k}, \eqref{c4k}, \eqref{c2k}, \eqref{ae16} and \eqref{ae15**} that
\begin{align}
b^k(\omega^k-\omega_{\Gamma}) &\rightharpoonup b(\omega-\omega_{\Gamma}) &&\textrm{weakly in } L^{q}(0,T;W^{1,q}_{\Gamma}(\Omega)) \textrm{ for all } q\in \left[1, \frac{16}{11}\right), \label{c11k*pppa}\\
\frac{b^k\omega^k}{1+b^k} &\rightharpoonup \frac{b\omega}{b+1} &&\textrm{weakly in } L^{q}(0,T;W^{1,q}(\Omega)) \textrm{ for all } q\in \left[1, \frac{16}{11}\right). \label{c11k*pppb}
\end{align}
Note that it also follows from \eqref{finalap}, \eqref{sc1k} and Fatou's lemma that
\begin{equation}
\sup_{t\in (0,T)} \left( \|b(t)\|_1 + \|\ln b(t)\|_1 \right) \le C. \label{bbln}
\end{equation}
Consequently, we see that the limit objects $(\bv,b,\omega,p,\bs)$ satisfy \eqref{k1T}--\eqref{k3P**} except \eqref{k2Tx}. Our goal is to prove \eqref{k2Tx} and all the remaining relations in Theorem~\ref{TH1}.

Thus, in order to identify also the limit in the nonlinear terms depending on $\omega^k$, we establish the strong convergence of the sequence $\omega^k$. To do so, we define two $4$-vectors
\begin{align*}
\ba^k&:=(\omega^k, \omega^k\bv^k - \mu^k \nabla \omega^k),\\
\bc^k&:=(b^k (1+b^k)^{-1} \omega^k, 0, 0, 0).
\end{align*}
Since, $\mu^k\le C b^k$ (due to the minimum principle for $\omega^k$), we use \eqref{minmaxok}, \eqref{ae14} and \eqref{ae2} to conclude that
$$
\|\ba^k\|_{L^{\frac{16-\lambda}{11}}(Q)} + \|\bc^k\|_{L^{\infty}(Q)} \le C(\lambda^{-1}) \qquad \textrm{ for all } \lambda \in (0,1).
$$
It also follows from \eqref{weak2-1k},  \eqref{minmaxok} and \eqref{ae16} that
\begin{align}
&\|\Div_{t,x} \ba^k\|_{L^{\infty}(Q)}=\|(\omega^k)^2\|_{L^{\infty}(Q)} \le C,\label{jidlo}\\
&\|\nabla_{t,x} \bc^k - (\nabla_{t,x} \bc^k)^T\|_{L^1(Q)}\le C\left\|\nabla_{x} \left(\frac{b^k\omega^k}{1+b^k}\right)\right\|_{L^1(Q)}\le C.\label{piti}
\end{align}
Then, it follows from \eqref{ae14} that
\begin{align}
\mu^k \nabla \omega^k &\rightharpoonup \overline{\mu \nabla \omega} &&\textrm{weakly in } L^{q}(0,T;L^{q}(\Omega)^3) \textrm{ for all } q\in \left[1, \frac{16}{11} \right) \label{c11k*p***}
\end{align}
and using the convergence results \eqref{sc1}, \eqref{sc3k} and \eqref{c2k} it is not difficult to observe that
\begin{align}
\ba^k &\rightharpoonup \ba &&\textrm{weakly in  } L^{q}(Q)^4 \textrm{ for all } q\in \left[1,\frac{16}{11}\right),\label{bak}\\
\bc^k &\rightharpoonup^* \bc &&\textrm{weakly$^*$ in } L^{\infty}(Q)^4, \label{bc1*kc}
\end{align}
where
\begin{align*}
\ba&:=(\omega, \omega\bv - \overline{\mu \nabla b}) \quad \textrm{ and } \quad \bc:=(b (1+b)^{-1} \omega, 0, 0, 0).
\end{align*}
Hence, we apply the Div--Curl Lemma (see \cite{Mu78, Ta78, Ta79}, \cite{feireisl04, FeNo09}) to the vector fields $\ba^k$ and $\bc^k$ and with the help of  \eqref{jidlo}--\eqref{bc1*kc} and also \eqref{c2k} we get
\begin{align}
\frac{b^k |\omega^k|^2}{1+b^k} &\rightharpoonup \frac{b |\omega|^2}{1+b} &&\textrm{weakly in  } L^{1}(Q).\label{somk}
\end{align}
This however directly implies (by using \eqref{sc1k}) that
\begin{align}
b^k \omega^k &\to b \omega &&\textrm{strongly in  } L^{2}(Q).\label{somks}
\end{align}
Since \eqref{bbln} implies that $b>0$ almost everywhere in $Q$, we finally deduce from \eqref{somks}, \eqref{minmaxok} and \eqref{sc1k} that
\begin{align}
\omega^k  &\to  \omega &&\textrm{strongly in  } L^{q}(0,T; L^q(\Omega)) \textrm{ for all } q\in [1,\infty).\label{stro}
\end{align}
Having \eqref{stro} we identify the limit of $\mu^k$: it follows from \eqref{sc1k} and the minimum principle \eqref{minmaxok} that
\begin{align}
\mu^k  &\to  \mu:=\frac{b}{\omega} &&\textrm{strongly in  } L^{q}(0,T; L^q(\Omega)) \textrm{ for all } q\in \left[1,\frac{8}{3}\right).\label{mustro}
\end{align}
As a direct consequence of \eqref{mustro}, \eqref{c1k} and \eqref{c4k} we then conclude that
\begin{align}
\mu^k \nabla b^k &\rightharpoonup \mu \nabla b &&\textrm{weakly in } L^{q}(0,T;L^{q}(\Omega)^3) \textrm{ for all } q\in \left[1, \frac87 \right), \label{c11k*p}\\
T_k(\mu^k) \bD(\bv^k) &\rightharpoonup \mu \bD(\bv) &&\textrm{weakly in } L^{q}(0,T;L^{q}(\Omega)^{3\times 3}) \textrm{ for all } q\in \left[1, \frac{16}{11}\right). \label{c11k*pp}
\end{align}
Finally, we focus on the identification of the weak limit in \eqref{c11k*p***}. First, it follows from \eqref{ae15} and \eqref{somks} and the uniqueness of the weak limit  that
\begin{align}
\nabla (b^k \omega^k) &\rightharpoonup \nabla (b \omega) \quad\textrm{weakly in } L^{q}(0,T;L^{q}(\Omega)^{3})\textrm{ for all }q\in \left[1, \frac{16}{11}\right). \label{ssopl}
\end{align}
Consequently, using the following identity (which is valid on the level of $k$-approximation valid, since $\omega^k$ is a Sobolev function)
$$
\mu^k \nabla \omega^k=\frac{\nabla (b^k\omega^k)}{\omega^k} - \nabla b^k,
$$
we observe that \eqref{minmaxok}, \eqref{c4k}, \eqref{c11k*p***},\eqref{stro} and \eqref{ssopl} finally imply that
\begin{align}
\mu^k \nabla \omega^k &\rightharpoonup \frac{\nabla (b\omega)}{\omega} - \nabla b &&\textrm{weakly in } L^{q}(0,T;L^{q}(\Omega)^3) \textrm{ for all } q\in \left[1, \frac{16}{11} \right). \label{dospol}
\end{align}
Moreover, by noting \eqref{ssopl} and \eqref{somks}, we can deduce similarly as above that
\begin{align}
b^k\omega^k &\to b\omega &&\textrm{strongly in } L^{1}(0,T;L^{1}(\partial \Omega)). \label{boundarydospol}
\end{align}
Using the assumptions \eqref{compat}, \eqref{lastass} and the definition of $\bg^k$ (see \eqref{dfgk} and combining them with the strong convergence results \eqref{sc1ktr} and \eqref{sc3ktr}, we see that $\bs$ fulfills \eqref{bc2} almost everywhere on $(0,T)\times \partial \Omega$. In addition, we get
\begin{equation}
\bg^k\cdot \bv^k \to -\bs \cdot \bv \quad \textrm{ a.e. on } (0,T)\times \Omega. \label{boundter}
\end{equation}
Thus, by using all of the above convergence results, we can easily let $k\to \infty$ in \eqref{weak1-1k} and \eqref{weak2-1k} to obtain \eqref{weak1-1kT} and \eqref{weak2-1kT}. Moreover, using the weak lower semicontinuity of norms, we can also let $k\to \infty$ in \eqref{weak3-1k} to obtain \eqref{entropy1}.

Next, we focus on establishing \eqref{weak3-1kT}. Denoting $E^k:=|\bv^k|^2/2 + b^k$, setting $\bw:=\bv^k z$ in \eqref{prk} with arbitrary $z\in W^{1,\infty}_{\Gamma}(\Omega)$ (which is now an admissible test function) and adding the result to \eqref{weak3-1k} we obtain
\begin{align}\label{ekeq}
&\left. \begin{aligned}
&\langle \partial_{t} E^k, z \rangle-((E^k+p^k) \bv^k, \nabla z) + (\bg^k(\cdot,b^k,\omega^k,\bv^k) \cdot \bv^k, z)_{\partial \Omega} \\
&\quad  + (\mu^k \nabla b^k, \nabla z)+\left(T_k(\mu^k)\bD(\bv^k)\bv^k, \nabla z\right)\\
&\quad =(-b^k \omega^k , z) +\frac12\left((2G_k(|\bv^k|^2)|\bv^k|^2-|\bv^k|^2 -\Gamma_k( |\bv^k|^2))\bv^k, \nabla z\right).
\end{aligned}\right.
\end{align}
Our goal is to study the limit of \eqref{ekeq} letting $k\to \infty$. First notice that it follows from \eqref{c3k*}, \eqref{sc3k} and the facts that $G_k(s)\to 1$ and $\Gamma_k(s)\to s$ as $k\to \infty$ that
\begin{align}
(2G_k(|\bv^k|^2)|\bv^k|^2-|\bv^k|^2 -\Gamma_k( |\bv^k|^2))\bv^k &\rightharpoonup \b0 &&\textrm{weakly in } L^{\frac{10}{9}}(0,T;L^{\frac{10}{9}}(\Omega)^3). \label{unava}
\end{align}
In the same manner, by using \eqref{sc3k} and \eqref{sc1k} (or \eqref{c4k*}), we conclude that
\begin{align}
\bv^kE^k &\rightharpoonup \bv E &&\textrm{weakly in } L^{\frac{10}{9}}(0,T;L^{\frac{10}{9}}(\Omega)^3). \label{unava1}
\end{align}
In addition, using \eqref{c11k*pp}, \eqref{c6kp} and \eqref{c3k*p}, we observe that
\begin{align}
T_k(\mu^k) \bD(\bv^k)\bv^k &\rightharpoonup \mu \bD(\bv)\bv \quad\textrm{weakly in } L^{q}(0,T;L^{q}(\Omega)^{3})\textrm{ for all }q\in \left[1, \frac{80}{79}\right), \label{c11k*ppo}\\
p_1^k \bv^k &\rightharpoonup p_1 \bv \quad\textrm{weakly in } L^{q}(0,T;L^{q}(\Omega)^{3})\textrm{ for all }q\in \left[1, \frac{80}{79}\right), \label{unava3}\\
p_2^k \bv^k &\rightharpoonup p_2 \bv \quad\textrm{weakly in } L^{q}(0,T;L^{q}(\Omega)^{3})\textrm{ for all }q\in \left[1, \frac{10}{9}\right). \label{unava4}
\end{align}
For the last term appearing in \eqref{ekeq} for which we have not identified the weak limit yet, we use \eqref{c4k*p*}, \eqref{c1k} and \eqref{sc3k} to conclude that
\begin{align}
\p_3^k bv^k  &\rightharpoonup p_3 \bv \quad\textrm{weakly in } L^{\beta_{\bg}}(0,T;L^{2}_{loc}(\Omega)^{3}). \label{unava5}
\end{align}
Hence, recalling also \eqref{c11k*p} and \eqref{somks}, we can deduce from \eqref{ekeq} that, for arbitrary $\Omega'\subset \overline{\Omega'} \subset \Omega$,
\begin{equation}
\int_0^T \|\partial_t E^k\|^q_{W^{-1,q}_0(\Omega')}\; dt \le C(q,\Omega') \qquad \textrm{ for all } q\in [1,\beta_{\min}), \label{et}
\end{equation}
where $\beta_{\min}$ is defined in Theorem~\ref{TH1}. Thus, taking $z\in \mathcal{D}(\Omega)$, we can let $k\to \infty$ in \eqref{ekeq} and obtain \eqref{weak3-1kT} (note that due to compact support of the test functions the boundary term vanishes and we can use better interior spatial regularity of $p_3$) and also \eqref{k2Tx}.

Finally, we focus on the global estimates for $E^k$. By using H\"{o}lder's inequality and the embedding theorem, we get (also with the help of the a~priori estimate \eqref{ae1})
\begin{equation}
\int_0^T\|\bv^k\|_{\frac{12}{5}}^8\; dt \le C\int_0^T \|\bv^k\|_{2}^6\|\bv^k\|_{6}^2\; dt\le C.\label{bizsi}
\end{equation}
Thus, assuming that $\beta_{\bg}>\frac87$, we have that
$$
\frac{1}{\beta_{\bg}}+\frac{1}{8}>1 \qquad \frac{2}{3\beta_{\bg}} + \frac{5}{12}>1.
$$
Consequently, it follows from \eqref{c4k*p}, \eqref{sc3k} and \eqref{bizsi} that there exists some $\beta_{1}>1$ such that
\begin{align}
\bv^k p_3^k &\rightharpoonup \bv p_3 \quad\textrm{weakly in } L^{\beta_1}(0,T;L^{\beta_1}(\Omega)^{3}). \label{unava8}
\end{align}
Thus, going back to \eqref{ekeq}, we can deduce that (here, if necessary, we possibly redefine $\beta_1$ so that $\beta_1<\beta_{\min}$)
\begin{equation}
\int_0^T \|\partial_t E^k\|^{\beta_1}_{W^{-1,\beta_1}_0(\Omega)}\; dt \le C. \label{www}
\end{equation}
Hence, after taking the limit $k\to \infty$ we can recover the relation \eqref{bettere}. In addition, we see that \eqref{weak3-1kT} holds for all $z\in L^{\infty}(0,T; W^{1,\infty}_0(\Omega))$. In addition, using the assumption \eqref{lowg} and the pointwise convergence \eqref{boundter}, we can use Fatou's lemma to derive \eqref{weak3-1kTE}.

Finally, if $\beta_{\bg}> \frac85$ then it follows from \eqref{ae2}, \eqref{ae3} that there exists a $\beta_2>0$ such that
$$
\int_0^T \int{\partial \Omega} |\bg^k\cdot \bv^k|^{\beta_2}\; dS \; dt\le C.
$$

Hence, we can improve the estimate \eqref{www} and we see that there is some $\beta_2>1$ such that
\begin{equation}
\int_0^T \|\partial_t E^k\|^{\beta_2}_{W^{-1,\beta_2}_{\Gamma}(\Omega)}\; dt \le C. \label{www2}
\end{equation}
Thus, the above relation allows us to show \eqref{bettere2} and, moreover, to prove the validity of \eqref{weak3-1kTE} with the equality sign. Attainment of the boundary conditions for $\bv_0$ and $\omega_0$ is classical. On the other hand, for $b_0$ we proceed differently. Here, we give only a very brief sketch and refer the interested reader to \cite{BuLeMa11} or \cite{BuMaRa09} for details. It can be shown with the help of \eqref{entropy1} that
$$
\sqrt{b(t)} \rightharpoonup \overline{\sqrt{b_0}} \ge \sqrt{b_0} \quad \textrm{ weakly in } L^2(\Omega).
$$
Moreover, it follows from \eqref{weak3-1kT} and from the fact that $\bv_0$ is attained that
$$
\lim_{t\to 0_+}\int_{\Omega}b(t)\le \int_{\Omega}b_0.
$$
From these two relationships one can deduce \eqref{oinda333T} for $b$. The proof of Theorem \ref{TH1} is thus complete.

%\begin{appendix}
%\section{Auxiliary results}
%\end{appendix}

%\bibliography{lipschitz}
%\bibliographystyle{plain}

\def\ocirc#1{\ifmmode\setbox0=\hbox{$#1$}\dimen0=\ht0 \advance\dimen0
  by1pt\rlap{\hbox to\wd0{\hss\raise\dimen0
  \hbox{\hskip.2em$\scriptscriptstyle\circ$}\hss}}#1\else {\accent"17 #1}\fi}
  \def\cprime{$'$}

\end{document}